\documentclass [12pt]{amsart}
\usepackage{latexsym,amssymb, amsmath}

\def\pf{\noindent\emph{Proof: }}
\def\stop{\hfill$\Box$}

\usepackage{amsmath, amsfonts, amssymb, amsthm}
\newtheorem{thm}{Theorem}

\newtheorem{lemma}[thm]{Lemma}

\newtheorem{defn}[thm]{Definition}
\newtheorem{prop}[thm]{Proposition}

\numberwithin{thm}{section}

\DeclareMathOperator{\Vol}{Vol}

\DeclareMathOperator{\Ric}{Ric}
\DeclareMathOperator{\Hess}{Hess}
\DeclareMathOperator{\diver}{div}

\begin{document}

\title [Compactness Results] {Some Compactness Results Related to Scalar
Curvature Deformation}

\author{Yu Yan}
\address{Department of Mathematics\\
         The University of British Columbia\\
         Vancouver, B.C., V6T 1Z2\\ Canada}
   \email{yyan@math.ubc.ca}

\begin{abstract}
Motivated by the prescribing scalar curvature problem, we study the equation $\Delta _g u +Ku^p=0 \,\ (1+\zeta
\leq p \leq \frac{n+2}{n-2})$ on locally
conformally flat manifolds $(M,g)$ with $R(g)=0$. We prove that when
$K$ satisfies certain conditions and the dimension of $M$ is 3 or 4, any solution $u$ of this equation
with bounded energy has uniform upper and lower bounds. Similar
techniques can also be applied to prove that on 4-dimensional scalar positive
manifolds the solutions of $\Delta _gu-\frac{n-2}{4(n-1)}R(g)u+Ku^p=0, K>0, 1+\zeta \leq
p \leq \frac{n+2}{n-2}$ can only have simple blow-up points.

\end{abstract}

\maketitle

\section{Introduction}

Let $(M^n,g)$ be an n-dimensional compact manifold with metric $g$,
and $u>0$ be a positive function defined on $M$. The scalar curvature of the
conformally deformed metric $u^{\frac{4}{n-2}}g$ is given by
$$R(u^{\frac{4}{n-2}}g)=-c(n)^{-1}u^{-\frac{n+2}{n-2}} \big (\Delta
_{g} u-c(n)R(g)u \big ) \hspace{.15in} \text{where }
c(n)=\frac{n-2}{4(n-1)}.$$

The famous Yamabe conjecture says that given a compact Riemannian
manifold $(M, g)$ of dimension $n\geq3$, $g$ can be conformally
deformed to a metric of constant scalar curvature. This conjecture was
proved to be true by the work of Trudinger (\cite{Tr}), Aubin (\cite{Au}) and
Schoen (\cite{S1}). 

It is natural to ask for a prescribed smooth function $K$ on $M$ if
it is possible to deform $g$ to a metric with scalar curvature
$K$. J. Escobar and R. Schoen studied this question in \cite{ES} and
gave some conditions under which $K$ can be a scalar curvature function.

Since the proofs of their results are variational, it is interesting
to know if the general compactness theorems hold under the same
conditions as well. On an n-dimensional
compact manifold $(M,g)$, are all the positive solutions of the
following equation compact in the $C^2$ norm? This
equation is the subcritical scalar curvature deformation equation 

$$
\Delta_{g}u - \frac{n-2}{4(n-1)} R(g)u +Ku^p=0 \hspace{.3in} \text{ where } 1<p\leq\frac{n+2}{n-2}
$$

If we can establish a uniform upper bound for the $C^0$ norm of $u$,
then the compactness in the $C^2$ norm will follow easily from the
bootstrap argument.

R. Schoen (\cite{S2}), and Y. Li and M. Zhu (\cite{LZ}) gained
compactness results on three dimensional manifolds not conformally
diffeomorphic to $S^3$ when $K$ is a positive constant and positive
function, respectively. 

In this paper we investigate the zero scalar curvature case.
This case is unknown and technically more difficult
than the positive scalar curvature case. When the scalar curvature is
zero, the equation becomes 
\begin{equation}
\label{eq:main}
\Delta _g u + Ku^p=0 \hspace{.3in} \text{ where } 1+\zeta \leq p \leq
\frac{n+2}{n-2}
\end{equation}
The necessary conditions for (\ref{eq:main}) to have solutions are $K>0$
somewhere on $M$ and $\int_{M} K dv_g <0$, hence $K$ has to change
signs on the manifold. The blow-up estimates used to prove the $K>0$
case depends on the lower bound on $K$, so it cannot be used on the
region where $K$ is small. 

One way to overcome this problem is to assume an energy bound on
$u$. Under this assumption, it can be proved that the maximum of $u$ is uniformly
bounded where $K \leq \delta$ for $\delta >0$ appropriately
small. This implies that blow-up can not happen and the $C^0$ norms
are bounded, which gives the compactness.

Furthermore, if the integrals $|\int_M K|$ have a uniform positive lower bound,
it can be proved that any limit function of a convergence sequence of positive
solutions of (\ref{eq:main}) is non-trivial, i.e. is strictly positive.

In summary, the main result in this paper is:

\begin{thm}
\label{thm:main}
Let $(M,g)$ be a three or four dimensional locally conformally flat compact
manifold with $R(g)=0$. Let $\mathcal{K}:=\{ K: K>0 $
somewhere on $M, \int_M K dv_g\leq -{C_K}^{-1}<0,$ and $  \| K
\|_{C^3}\leq C_K \}$ for some constant $C_K$, and
$S_{\Lambda} :=\{u: u>0 $ solves $(\ref{eq:main})$, $K \in \mathcal{K}, $ and $ E(u)
\leq \Lambda \}.$ Then there exists $C=C(M, g, C_K, \Lambda, \zeta)>0$ such
that $u \in S_{\Lambda}$ satisfies $C^{-1} \leq \|u\|_{C^3(M)}\leq C.$

\end{thm}

\noindent
This theorem is consistent with the existence theorem of Escobar
and Schoen. 

Additionally, we can use some of the techniques in the proof of the
above theorem to get a better understanding of the possible blow-up for 4
dimensional scalar positive manifolds.

\begin{thm}
\label{thm:cor}
If $(M,g)$ is a four dimensional locally conformally flat compact
manifold with positive scalar curvature, 
and $K>0$ on $M$, then the possible blow-up of the solutions
of equation 
\begin{equation}
\Delta _g u- c(n)R(g)u + Ku^p=0 \hspace{.3in} (1 \leq p\leq \frac{n+2}{n-2})
\end{equation}
 is always simple.
\end{thm}

\noindent
We will give the precise definition of simple blow-up in section
5. Roughly it means that blow-up points are isolated and consist of a
simple ``bubble''. 

The rest of this paper is mostly devoted to the proof of Theorem \ref{thm:main}, and
it will also illustrate the proof of Theorem \ref{thm:cor} .
In section \ref{section:lower} we prove the lower bound on $u$ assuming the upper bound
exists. In section \ref{section:smallK} we show that $u$ is uniformly bounded above on
the region where $K$ is sufficiently small. In section
\ref{section:iso} we reduce the
possible blow-up on the region where $K$ is big to two cases. In
section \ref{section:simple} we introduce the definition of simple blow-up and prove some
important estimates. In sections \ref{section:case1} and \ref{section:case2} we show that neither case in
section \ref{section:iso} can happen, hence prove the compactness theorem.

\section{The Lower Bound on $u$}
\label{section:lower}

Suppose $u$ with $E(u)\leq \Lambda$ is a solution of equation (\ref{eq:main}) for some $1+\zeta \leq p
\leq \frac{n+2}{n-2}$, where $K$ satisfies $$ K>0 \text { somewhere on
} M, \hspace{.1in} \|K\|_{C^3}\leq C_K, \hspace{.1in}
\text{and} \hspace{.1in} \int_M K dv_g \leq -\frac{1}{C_K}<0.$$ 

\begin{lemma}
\label{lem:lb}
$\int _M u^2 dv_g \leq C \int _M |\triangledown u|^2 dv_g $
where $C=C(M, g, n, C_K).$
\end{lemma}

\pf
Let $\bar{u}=\Vol(M)^{-1}\int _{M} u \,\,dv_g$, then

\begin{equation}
\label{eq:lbsep}
\int _M u^2 dv_g 
 \leq  2 \int _M (u-\bar{u})^2 dv_g + 2 \int _M \bar{u}^2 dv_g
\end{equation}
By the Poincar\'{e} inequality
\begin{equation}
\label{eq:Poin}
\int _M (u-\bar{u})^2 dv_g \leq C(M,g) \int _M |\triangledown u|^2 dv_g, 
\end{equation}
so we only need to find an upper bound for $\bar{u}$. 

\begin{eqnarray*}
& & \bar{u}^{p+1}C_K^{-1}\\
 &\leq & \bar{u}^{p+1}(\int _M -K dv_g)\\
&=& -\int _M Ku^{p+1}dv_g + \int _M
K(u^{p+1}-\bar{u}^{p+1})dv_g \\
& \leq & \max _M |K| \int _M \big | u^{p+1}-\bar{u}^{p+1} \big |dv_g\\
&\leq & \max _M |K| \int _M \left (2^{p-1}(p+1)|u-\bar{u}|^{p+1}+2^{2p-1}(p+1)\bar{u}^p|u-\bar{u}| \right)dv_g\\
& & \text{(by calculus)}\\
&\leq & C(n, C_K) \int _M \left( |u-\bar{u}|^{p+1} + \bar{u}^p |u-\bar{u}|  \right )dv_g.
\end{eqnarray*}
Therefore $$\bar{u}^{p+1} \leq C(n, C_K) \int _M
\left ( |u-\bar{u}|^{p+1}+\bar{u}^p|u-\bar{u}| \right )dv_g.$$

\noindent
The first term on the right hand side
$$
\int _M |u-\bar{u}|^{p+1}dv_g \leq C(M, g, n) \left (\int _M |\triangledown u|^2 dv_g \right )^{\frac{p+1}{2}}
$$
by H\"{o}lder and Sobolev inequalities.

\noindent
Similarly the second term on the right hand side
$$
\int _M |u-\bar{u}| dv_g \leq  C(M, g, n) \left (\int_M |\triangledown u|^2 dv_g \right )^{\frac{1}{2}}.
$$

\noindent
Therefore $$ \bar{u}^{p+1} \leq  C(M, g, n, C_K)\left ( \left (\int_M |\triangledown u|^2 dv_g \right )^{\frac{p+1}{2}}+\bar{u}^p
\left (\int_M |\triangledown u|^2 dv_g \right )^{\frac{1}{2}} \right) $$
Choose $C=\max\{1, C(M, g, n, C_K)\},$ then $\bar{u} \leq C
\left (\int_M |\triangledown u|^2 dv_g \right ) ^{\frac{1}{2}}$. Thus
by (\ref{eq:lbsep}) and (\ref{eq:Poin})
$$\int _M u^2 dv_g \leq C(M, g, n, C_K) \int_M |\triangledown u|^2 dv_g .$$
 
\stop \\

Now assume $u$ is bounded above, we claim that it
is also bounded below away from 0. Suppose not, then $\exists \,\, \{x_i\} \subset M, \{ u_i \} \subset S_{\Lambda}$ and $\{ K_i
\} \subset \mathcal{K}$ such that $\Delta u_i + K_i u_i ^{p_i}=0$ and
$u_i(x_i) \rightarrow  0$. Since $\|u _i\|_{C^3}$ is
bounded, then there is a subsequence also denoted as $\{ u_i \}$ which
converges in $C^2$-norm to some function $u \geq 0$. Similarly after
passing to a subsequence $\{ K _i \} $ also converges to
some function $K$. Let $p$ be the corresponding limit point of $\{ p_i
\} $, then we have $$\Delta u + Ku^p=0.$$ 
The manifold $M$ is compact,
so $\{x_i\}$ also has a limit point $x \in M$. Since $u_i(x_i) \rightarrow
0$, $u(x)=0$ and then by
the strong maximum principle $u \equiv 0$.

\noindent
On the other hand, by the Sobolev inequality and Lemma \ref{lem:lb}
\begin{eqnarray*}
\left (\int _M u_i^{\frac{2n}{n-2}} dv_g \right )^{\frac{n-2}{n}} 
 & \leq  & C\int_M |\triangledown u_i|^2 dv_g \\
& = & C \int _M K_i u_i^{p_i+1}dv_g \\
& \leq & C \left (\int _M u_i ^{\frac{2n}{n-2}}dv_g \right )^{\frac{n-2}{2n} (p_i+1)}
\end{eqnarray*}
This implies $$1 \leq C \left (\int _M
u_i^{\frac{2n}{n-2}}dv_g \right )^{\frac{n-2}{n} \frac{p_i-1}{2}}.$$  Then since $p_i-1 \geq \zeta>0$, $u_i$ cannot converge to 0, contradicting $u \equiv 0$. 

Now it is only left to show that $ u $ has
a uniform upper bound. Since $u$ satisfies equation (\ref{eq:main}), the upper
bound of $\| u \| _{C^3}$ will follow easily by the standard
elliptic theory and Sobolev embedding theorem once we establish a
uniform upper bound on $u$.

\section{An Upper Bound on $u$ on the set where $K$ is Small}
\label{section:smallK}

Let $u>0$ be a function which satisfies 
$$ \Delta u + Ku^p =0, \hspace{.3in} 1 \leq p \leq
\frac{n+2}{n-2}$$ for some $K \in \mathcal{K}$. 

Choose $1 < \beta < \frac{2n}{n-2} -1$. Let $x$ be a point on $M$ and $\varphi$ be a cut-off
function such that $$\varphi \equiv 1 \text{ on } B_{\frac{3}{4}
\sigma}(x), \hspace{.1in} \varphi \equiv 0 \text{ on } M
\setminus B_{\sigma}(x), \hspace{.1in} \text{and } \,\,
|\triangledown \varphi | \leq 2 \sigma ^{-2}.$$ 

\noindent
Multiplying (\ref{eq:main}) by $\varphi ^2 u^{\beta}$ and integrating by parts
gives

\begin{eqnarray*}
& & \int_M \varphi ^2 K u^{\beta +p}dv_g\\
 &=&-\int_M \varphi ^2 u^{\beta}\Delta _g u \,\,dv_g\\
&=&\frac{4\beta}{(\beta+1)^2} \int_M \varphi ^2 |\triangledown (u
^{\frac{\beta +1}{2}})|^2  dv_g 
 + \frac{2}{\beta +1} \int _M
u^{\frac{\beta +1}{2}} \triangledown (u^{\frac{\beta +1}{2}})\cdot
\triangledown (\varphi ^2) dv_g
\end{eqnarray*}

\noindent
Let $w=u^{\frac{\beta+1}{2}}$, then 
\begin{eqnarray*}
& & \frac{4\beta}{(\beta +1)^2} \int _M \varphi ^2  |\triangledown w|^2 \,\,
dv_g\\
 & = & -\frac{4}{\beta +1} \int _M \varphi w \triangledown w \cdot
\triangledown \varphi \,\, dv_g + \int _M \varphi^2 K w^2 u^{p
-1} \,\, dv_g\\
&\leq & \frac{4}{\beta +1}\epsilon \int_M \varphi ^2 |\triangledown w|^2\,\, dv_g +
\frac{4}{\beta+1}\epsilon ^{-1} \int_M w^2 |\triangledown \varphi|^2\,\, dv_g \\
& & + \max _{B_{\sigma}(x)}K\int _M \varphi ^2 w^2 u^{p-1} \,\, dv_g 
\end{eqnarray*}

\noindent
Choosing $\epsilon$ small enough (only depending on $\beta$) we can
absorb the first integral into the left hand side and get
$$\int_M \varphi ^2 |\triangledown w|^2 \,\, dv_g \leq C \sigma ^{-2}
\int _{B_{\sigma}(x)} w^2 \,\, dv_g +  \max _{B_{\sigma}(x)}K \int _M \varphi ^2 w^2 u^{p-1}
\,\, dv_g.$$
Therefore 
\begin{eqnarray*}
\int_M |\triangledown (\varphi w)|^2 dv_g &\leq & 2 \int_M
|\triangledown \varphi|^2 w^2 \,\, dv_g + 2\int _M \varphi ^2
|\triangledown w|^2 dv_g \\
&\leq& C \sigma ^{-2} \int _{B_{\sigma}(x)} w^2 dv_g + 2 \max _{B_{\sigma}(x)}K \int _M
\varphi ^2 w^2 u^{p-1} dv_g
\end{eqnarray*}
Since $\beta+1<\frac{2n}{n-2}$,
\begin{equation}
\label{eq:w2bd}
\int_{B_{\sigma}(x)} w^2 dv_g \,\, \leq  \,\,  \int_{B_{\sigma}(x)}
u^{\frac{2n}{n-2}}\,\, dv_g + \Vol(M,g) \,\, \leq \,\, C(M, g, n, C_K, \Lambda)
\end{equation}
by the Sobolev inequality, Lemma \ref{lem:lb} and the energy bound on $u$. 

\noindent
Similarly
$(p-1)\frac{n}{2} \leq \frac{2n}{n-2}$ implies 

\begin{eqnarray*}
& & \int _M \varphi ^2 w^2 u^{p-1}\,\,dv_g\\
 &\leq& \left (\int _{B_{\sigma}(x)} (\varphi ^2
w^2)^{\frac{n}{n-2}}\,\,dv_g \right )^{\frac{n-2}{n}} \left (\int_{B_{\sigma}(x)}
(u^{p-1})^{\frac{n}{2}} \,\, dv_g \right )^{\frac{2}{n}}\\
& \leq & C(M, g, n, C_K, \Lambda) \left (\int _{B_{\sigma}(x)} (\varphi
w)^{\frac{2n}{n-2}}\,\, dv_g \right )^{\frac{n-2}{n}}.
\end{eqnarray*} 
So we know

\begin{eqnarray}
\label{eq:gradw}
\int_M |\triangledown (\varphi w)|^2 dv_g 
& \leq & C(M,g,n,C_K,
\Lambda)\bigg [ \sigma ^{-2}  +  \nonumber \\
& &  \max _{B_{\sigma}(x)}K \left (\int _{B_{\sigma}(x)} (\varphi
w)^{\frac{2n}{n-2}}\,\, dv_g \right )^{\frac{n-2}{n}}  \bigg ]
\end{eqnarray}

\noindent
Then by the Sobolev inequality
\begin{eqnarray*}
& & \left (\int _M (\varphi w)^{\frac{2n}{n-2}}\,\, dv_g \right )^{\frac{n-2}{n}}\\
 &\leq &
C(M, g) \left (\int _M |\triangledown (\varphi w)|^2 \,\, dv_g + \int _M
(\varphi w)^2 \,\,dv_g \right )\\ 
&\leq& C(M,g, n, C_K, \Lambda)\bigg[ \sigma ^{-2} +  \max _{B_{\sigma}(x)}K \left (\int _{B_{\sigma}(x)} (\varphi
w)^{\frac{2n}{n-2}}\,\, dv_g \right )^{\frac{n-2}{n}} \bigg ]  
\end{eqnarray*}
where the last inequality follows from (\ref{eq:w2bd}) and (\ref{eq:gradw}).

Let $\delta = \frac{1}{4}C(M, g, n, C_K, \Lambda)^{-1}$. If $\max _{B_{\sigma}(x)}K < 2\delta$, then we can absorb the
second term on the right hand side of the above inequality into the
left hand side to get $$
\left ( \int _M (\varphi w)^{\frac{2n}{n-2}}\,\, dv_g \right )^{\frac{n-2}{n}} \leq
C(M,g,n,C_K, \Lambda) \sigma ^{-2}.$$
Therefore 
\begin{equation}
\label{eq:control-upper-u}
\int _{B_{\frac{3}{4}\sigma}(x)} u^{\frac{\beta +1}{2}\frac{2n}{n-2}}\,\, dv_g\leq
C(M,g,n,C_K, \Lambda) \sigma ^{-\frac{2n}{n-2}}
\end{equation}

Define $K_{\delta} := \{x \in M : K(x) <\delta\}$, and let $x_1 \in
\partial K_{\delta}$ and $x_2 \in \partial K_{2\delta}$ be the points which realize the
distance between $\partial K_{\delta}$ and $\partial K_{2\delta} $, i.e., $d _g(x_1, x_2) = d _g
(\partial K_{\delta}, \partial K_{2\delta}).$ Then 
\begin{eqnarray*}
2\delta = K(x_2)
& \leq & K(x_1)+ \max _M |\triangledown K| d_g(x_1, x_2)\\
&\leq& \delta + C_K d_g(x_1, x_2)
\end{eqnarray*}
which implies $$d_g(x_1, x_2)>C( \delta, C_K)=C(M, g, n, \Lambda, C_K).$$

\noindent
Let $\sigma = \frac{1}{2}d_g(x_1, x_2)>\frac{1}{2} C(M, g,n, \Lambda, C_K).$ For any $x \in K_{\delta}$,
we have $B_{\sigma}(x) \subset K_{2\delta}$, therefore
\begin{equation}
\label{eq:uLrbd}
\int _{B_{\frac{3}{4}\sigma}(x)} u^{\frac{\beta +1}{2}\frac{2n}{n-2}}\,\, dv_g\leq
C(M,g,n, \Lambda, C_K)
\end{equation}
by (\ref{eq:control-upper-u}) and the lower bound on $\sigma$.

\noindent
This tells us that $u \in L^{\frac{\beta +1}{2} \frac {2n}{n-2}
} \left (B_{\frac{3}{4}\sigma}(x) \right )$, hence $Ku^{p-1} \in
L^r \left (B_{\frac{3}{4}\sigma}(x) \right )$ where 
\begin{eqnarray}
\label{eq:rbd}
r \,\,= \,\, \frac{\beta +1}{2} \frac {2n}{n-2} \frac{1}{p-1} \,\,
\geq \,\, \frac{\beta+1}{2} \frac{n}{2} \,\, >  \,\, \frac{n}{2}&  \text{ when } & p\neq 1  \\ 
\text { and let } \hspace{.1in} r=\frac{\beta+1}{2}\cdot \frac{n}{2} & \text { when } & p=1. \nonumber
\end{eqnarray}

\noindent
By the elliptic theory 
\begin{equation}
\label{eq:supusmK}
\sup _{B_{\frac{1}{2}\sigma}(x)} u \leq C \left (M, g, r, 
\|Ku^{p-1}\|_{L^r \left ({B_{\frac{3}{4}\sigma}(x)} \right)}\right )\sigma ^{-\frac{n}{2}}
\|u\|_{L^2 \left (B_{\frac{3}{4}
\sigma}(x) \right )} 
\end{equation}

\noindent
The constant in the above inequality usually blows up when $r \rightarrow \frac{n}{2}$ or when $\|Ku^{p-1}\|_{L^r
\left ({B_{\frac{3}{4}\sigma}(x)} \right )}$ is unbounded. But here we have a
fixed lower bound on $r$ from (\ref{eq:rbd}), and by (\ref{eq:uLrbd})
$$
\| Ku^{p-1} \| _{L^r \left ({B_{\frac{3}{4}\sigma}(x)}\right )} \leq C(M, g, n, \zeta,
\Lambda, C_K). $$ 
So the constant in (\ref{eq:supusmK}) has an upper bound only depending on $M, $$ g, $$
n, $$ \zeta, \,\, $ $ \Lambda $ and $ C_K$.

\noindent
Since we also have a uniform lower bound on $\sigma$, and by Lemma \ref{lem:lb} and the energy bound on $u$ we know that
\begin{eqnarray*}  
\|u\|^2_{L^2 \left (B_{\frac{3}{4} 
\sigma}(x) \right )} 
&\leq& C(M, g, n, \zeta, C_K)\Lambda,
\end{eqnarray*}
it can be concluded from (\ref{eq:supusmK}) that  
$$\sup _{B_{\frac{1}{2}\sigma}(x)} u \leq C(M, g, n, \zeta, \Lambda, C_K). $$
Therefore
$$\sup _{K_{\delta}} u \leq C(M, g, n, \zeta, \Lambda, C_K) $$ since $x$
is an arbitrary point in $K_{\delta}$.

In the next few sections we are going to prove that $u$ is also
uniformly bounded on $M \setminus K_{\delta}$.

\section{Reduction to the Isolated Blow-Up Case}
\label{section:iso}

We first prove a lemma.

\begin{lemma}
\label{lem:iso}
Suppose $W \supseteq K_{\frac{\delta}{2}}$ ($\delta$ is chosen as in
section $\ref{section:smallK}$ and $K_{\frac{\delta}{2}}$ is also defined in the same way as in
section $\ref{section:smallK}$) is a compact subset of $M$ and $u>0$ is a solution of
equation {\em (\ref{eq:main})} with $K \in \mathcal{K}$ (as defined in
Theorem $\ref{thm:main} $). Given $\epsilon, R>0$, there
exists $C=C(\epsilon, R)>0$ such that if $\max _{\overline{M \setminus
W}} d(x)^{\frac{2}{p-1}}u(x) \geq C$, where $d(x)=dist _g(x, W)$
(let $d(x)=1$ in case $W=\varnothing$), then 
\begin{itemize}
\item $\frac{n+2}{n-2} -p <
\epsilon$ 
\item there exists $x_0 \in M \setminus W$ which is a
local maximum point of $u$, and the geodesic ball $B_{Ru(x_0)^{-\frac{p-1}{2}}}(x_0) \subset
\overline{M\setminus W}$ 
\item Choose the $y$-coordinates around $x_0$ so 
that $z=\frac{y}{u(x_0)^{\frac{p-1}{2}}}$ is a geodesic normal
coordinate system centered at $x_0$, then  $$ \Bigg \| u(0)^{-1}u \left (\frac
{y}{u(0)^{\frac{p-1}{2}}} \right ) -\bar{v}(y) \Bigg \|_{C^2(B_{2R}(0))} <
\epsilon$$
where $$\bar{v} (y)=\left (1+
\frac{K(x_0)}{n(n-2)}|y|^2 \right )^{-\frac{n-2}{2}}$$ is the (unique) exact solution
of $$\Delta v(y) + K(x_0)v(y)^{\frac{n+2}{n-2}}=0,\hspace{.1in} y\in
\mathbf {R}^n,\hspace{.1in} 0<v\leq 1, \hspace{.1in} v(0)=1 $$where $\Delta$
is the Euclidean Laplacian. 
\end{itemize}
\end{lemma}

\pf
Suppose no such $C$ exists, then there exists $\{ u_i\}, \{p_i\}, \{K_i\}$
and $\{ W_i \}$ such that $$\Delta _g u_i + K_i u_i ^{p_i}=0, $$$$ \max _{\overline{M \setminus W_i}} d_i(x)^{\frac{2}{p_i-1}}
u_i(x) \geq i \hspace{.1in}\text{ where } \hspace{.1in} d_i(x)=dist_g(x, W_i),$$ and
$$W_i \supseteq K_{\frac{\delta}{2}} ^{(i)} := \{x \in M : K_i(x) < \frac{\delta}{2} \},$$
but for each $i$ there doesn't exist any local maximum point which satisfies the conditions in
the lemma.

Define $f_i(x)=d_i(x)^{\frac{2}{p_i-1}}u_i(x)$, let $x_i$ be a maximum
point of $f_i(x)$ in $\overline{M \setminus W_i}.$ Since on the
boundary of $M \setminus W_i$, $f_i \equiv 0$, we know $x_i \in M
\setminus W_i$. Let $z$ be the geodesic normal coordinates with
respect to the background metric $g$ centered at $x_i$, and let
$y=u_i(x_i)^{\frac{p_i-1}{2}}z$. Define
$$v_i(y)=u_i(x_i)^{-1}u_i \left (\frac{y}{u_i(x_i)^{\frac{p_i-1}{2}}} \right ) .$$
Then $v_i$ satisfies 
\begin{equation}
\label{eq:rseq}
\Delta _{g^{(i)}} v_i + K_i\left (\frac{y}{u_i(x_i) ^{\frac{p_i-1}{2}}}
\right )v_i^{p_i}=0
\end{equation}
where the metric
$g^{(i)}(y)=g_{\alpha \beta} \left
(\frac{y}{u_i(x_i)^{\frac{p_i-1}{2}}}\right
)dy^{\alpha}dy^{\beta}$. Since we have chosen coordinates for which
$g_{\alpha \beta}(0)= \delta _{\alpha \beta}$, for
$y$ on a bounded set, $g^{(i)}(y)$ converges to the Euclidean metric.\\
Let $r_i =\frac{d_i(x_i)}{2}$ and
$R_i=u_i(x_i)^{\frac{p_i-1}{2}}r_i$, then by our choice of $x_i$, $R_i \geq \frac{1}{2}i^{\frac{p_i-1}{2}} \rightarrow \infty $ as $i \rightarrow \infty$.
In the ball centered at $x_i$ with radius $r_i$,
$$\frac{\max_{B_{r_i}}d_i(x)}{\min_{B_{r_i}}d_i(x)}=\frac{d_i(x_i)+\frac{d_i(x_i)}{2}}{d_i(x_i)-\frac{d_i(x_i)}{2}}=3.$$

\noindent
Also from the choice of $x_i$ we know$$ d_i(x)^{\frac{2}{p_i-1}}u_i(x)=f_i(x)
\leq f_i(x_i) = d_i(x_i)^{\frac{2}{p_i-1}}u_i(x_i),$$ so $$u_i(x)
\leq 3^{\frac{2}{p_i-1}}u_i(x_i), \hspace{.2in} \text{ i.e. } \hspace{.2in}
u_i \left (\frac{y}{u_i(0)^{\frac{p_i-1}{2}}} \right ) \leq 3^{\frac{2}{p_i -1}}u_i(0)$$
in the $y$ coordinates. Therefore 
\begin{equation}
\label{eq:ubvi}
v_i(y) \leq 3^{\frac{2}{p_i-1}}
\leq C \hspace{.2in} \text{ on } B_{R_i}(0) \subset \mathbf{R}^n(y) 
\end{equation}
for some constant $C$ independent of $i$.

By our choice of $x_i$, $d_i(x_i)^{\frac{2}{p_i-1}}u_i(x_i) \rightarrow
\infty$ as $i \rightarrow \infty$, but $d_i(x_i) $ is bounded above on
$M$, so $u_i(x_i) \rightarrow \infty$ as $ i \rightarrow
\infty$. Hence 
 $ \Big |\frac{y}{u_i(x_i)^{\frac{p_i-1}{2}}} \Big | \rightarrow 0  $ on any
compact subset on the $y$-plane.

\noindent
After passing to a subsequence, we may assume $x_i \rightarrow x_0 \in M $ and
$K_i$ converges to some function $K$ in $C^2$-norm. By the choice of $x_i$ and $W_i$ we know
$K_i(x_i) > \frac{\delta}{2}$, so $K(x_0) \geq \frac{\delta}{2}$.

\noindent
Then with (\ref{eq:rseq}) and (\ref{eq:ubvi}) by the elliptic
estimates we know $\{v_i\}$ converges to some function $\hat{v}$ in
$C^2$-norm which satisfies $\hat{v}(y) \geq 0$ and 
\begin{equation}
\Delta \hat{v}(y) + K(x_0) \hat{v}(y) ^{\frac{n+2}{n-2}} =0 \hspace{.1in}
\text{ on } \mathbf{R}^n(y)
\end{equation}   
Here and throughout the rest of this chapter we use $\Delta$ to denote
the Laplacian with respect to the Euclidean metric. The reason for
$\lim _{i \rightarrow \infty} p_i =\frac{n+2}{n-2}$ is that
$K(x_0)>0$, $v(0)=1$ (since $v_i(0) \equiv 1$) and there is no non-trivial solution of
$\Delta v+cv^p=0$ on $\mathbf{R}^n$ with $c>0$ and $1\leq p <
\frac{n+2}{n-2}$ (proved in \cite{CGS}).
Another consequence of  $v(0)=1$ is that by the Harnack inequality $v(y) >0$ for any $y$.

As proved in \cite{CGS}, $\hat{v}(y)$ has the expression $$
\hat{v}(y)=\left (\frac{\lambda
K(x_0)^{-\frac{1}{2}}}{1+\frac{1}{n(n-2)}\lambda ^2 |y-\bar{y}|^2
}\right)^{\frac{n-2}{2}} $$ for some $\lambda >0$, where $\bar{y}$ is the maximum point of $\hat{v}$.

\noindent
Since $v_i \rightarrow \hat{v}$ in $C^2$-norm on compact subsets on
$\mathbf{R}^n(y)$, there exists $\{ \bar{y}_i  \}$ such that each $\bar{y}_i$ is a local maximum
point of $v_i(y)$ and $\bar{y}_i \rightarrow \bar{y}$. Let $\bar{z} _i =u_i(x_i)^{-\frac {p_i-1}{2}}
\bar{y}_i $ and $\bar{x}_i =\exp _{x_i} \bar{z} _i \in M$. Note that
$R_i \rightarrow \infty$, so for large enough $i$, $\bar{y} _i \in B_{R_i}(0) \subset \mathbf{R}^n(y)$. Therefore $\bar{z} _i
\in B_{r_i}(0) \subset \mathbf{R}^n(z)$, hence $\bar{x} _i \in B_{r_i}(x_i)$, consequently $B_{r_i}(\bar{x}_i) \in M
\setminus W_i$ and $K_i(\bar{x} _i) > \frac{\delta}{2} $.

Let $l$ be any fixed large radius so that $|\bar{y}| < l$. Since $R_i =r_i u_i(x_i)^{\frac{p_i -1}{2}}$ $ \rightarrow \infty$, for
large enough $i$, $2l < R_i$. Hence $$
B_{\frac{2l}{u_i(x_i)^{\frac{p_i-1}{2}}}}(x _i) \subset B_{r_i}(x_i)
\subset M \setminus W_i.   $$

\noindent
When $i$ is large enough $\bar{y} _i $ is also in
$B_{l}(0)$, so $|\bar{z}_i| <
\frac{l}{u_i(x_i)^{\frac{p_i-1}{2}}}$. This implies that
$$B_{\frac{l}{u_i(x_i)^{\frac{p_i-1}{2}}}}(\bar{x} _i) \subset
B_{\frac{2l}{u_i(x_i)^{\frac{p_i-1}{2}}}}(x _i) \subset M
\setminus W_i.   $$

\noindent
We chose $\bar{y}_i$ to be a maximum point of $v_i(y)$, and
$\bar{y}_i \rightarrow \bar{y}$ which is the only maximum point of
$\hat{v}$. Thus when $i$ is sufficiently large $\bar{y}_i$ is also the
only maximum point of $v_i$ for $|y| \leq 2l$. Therefore $\bar{x}_i$ is the maximum
point of $u_i$ on $B_{\frac{2l}{u_i(x_i)^{\frac{p_i-1}{2}}}}(x_i)$.
In
particular $u_i(x_i) \leq u_i(\bar{x}_i)$, so $$
B_{\frac{l}{u_i(x_i)^{\frac{p_i-1}{2}}}}(\bar{x}_i) 
\supset B_{\frac{l}{u_i(\bar{x}_i)^{\frac{p_i-1}{2}}}}(\bar{x}_i).
 $$ Hence 
\begin{itemize}
\item $B_{\frac{l}{u_i(\bar{x}_i)^{\frac{p_i-1}{2}}}}(\bar{x}_i) \subset M
\setminus W_i$.
\end{itemize} 

Now redefine $z$ to be the geodesic normal coordinates centered at
each $\bar{x}_i$ and let $y=u_i(\bar{x}_i)^{\frac{p_i-1}{2}}z$. Define
$$\bar{v} _i(y) =
u_i(\bar{x}_i)^{-1}u_i \left
(\frac{y}{u_i(\bar{x}_i)^{\frac{p_i-1}{2}}} \right ).$$
Since $\bar{x}_i$ is the maximum point of $u_i$ on
$$B_{\frac{2l}{u_i(x_i)^{\frac{p_i-1}{2}}}}(x_i) \supset B_{\frac{l}{u_i(x_i)^{\frac{p_i-1}{2}}}}(\bar{x}_i) 
\supset B_{\frac{l}{u_i(\bar{x}_i)^{\frac{p_i-1}{2}}}}(\bar{x}_i),$$
then for all $|y| \leq l$, we have $\bar{v}_i(y) \leq 1$. 

\noindent
Then the same argument as that for $v_i$ shows that $\bar{v} _i$
converges in $C^2$-norm to some function $\bar{v}>0$ which satisfies 
\begin{equation}
\Delta \bar{v} + K(\bar{x}_0) \bar{v} ^{\frac{n+2}{n-2}}=0  
\end{equation}
where $\bar{x}_0=\lim _{i \rightarrow \infty} \bar{x} _{i}$ and
$K(\bar{x}_0) >0$ (because $K_i(\bar{x} _i) > \frac{\delta}{2}$).

\noindent
Since $y=0$ is a maximum point of $\bar{v}_i$, the maximum of
$\bar{v}$ is attained at $y=0$. So $\bar{v}$ has the expression
$$\bar{v}(y)= \left (1+\frac{K(\bar{x}_0)}{n(n-2)} |y|^2 \right )^{-\frac{n-2}{2}}. $$
Then since $\bar{v} _i \rightarrow \bar{v}$ and $\bar{x} _i
\rightarrow \bar{x} _0$, for large enough $i$,
\begin{itemize}
\item $
\bigg \|
u_i(\bar{x}_i)^{-1}u_i\left (\frac{y}{u_i(\bar{x}_i)^{\frac{p_i-1}{2}}}
\right ) - \left
(1+\frac{K_i(\bar{x}_i)}{n(n-2)} |y|^2 \right )^{-\frac{n-2}{2}} \bigg \| <
\epsilon.  $
\end{itemize}
This is a contradiction. 

\stop \\

Now fix $\epsilon >0$ and $R >>0$. Suppose $$\max _{\overline{M \setminus
K_{\frac{\delta}{2}}}} d_g \left (x, M \setminus
K_{\frac{\delta}{2}} \right )^{\frac{2}{p-1}}u(x) > C \hspace{.1in} (C \text{
is the constant in Lemma } \ref{lem:iso}). $$  
Applying Lemma \ref{lem:iso} to the case $W=K_{\frac{\delta}{2}}$, there exists $x_1 \in
M$ which is a local maximum point of $u$ and satisfies the conditions
in the lemma.\\
Let $r_1=Ru(x_1)^{-\frac{p-1}{2}}$. 
We can stop the procedure if at any point $x \in M$, $d_g \left (x, \overline{K_{\frac{\delta}{2}}
\bigcup B_{r_1}(x_1)} \right )^{\frac{2}{p-1}}u(x) $ is bounded by some
constant only depending on $\epsilon$ and $R$. Otherwise let $W=\overline{K_{\frac{\delta}{2}}
\bigcup B_{r_1}(x_1)}$  and apply
Lemma \ref{lem:iso} again to find another local maximum point $x_2 \in M$ to
satisfy the conditions in the lemma.

\noindent
Repeating this procedure we will get a sequence of disjoint balls
$$B_{r_1}(x_1), ..., B_{r_N}(x_N).$$ The
sequence must be finite 
because for each $i$ $$r_i =R u(x_i)^{-\frac{p-1}{2}} \geq R (\max _M
u)^{-\frac{p-1}{2}},$$ hence the volume of $B_{r_i}(x_i)$ has a lower
bound. Here we allow the number of balls $N$ to depend on the function $u$.\\
Thus we know $$u(x) \leq C d_g \left (x, \overline{K_{\frac{\delta}{2}}
\bigcup \cup _{i=1} ^{N}B_{r_i}(x_i)}\right )^{-\frac{2}{p-1}}$$ for some
constant $C=C(\epsilon, R)$.

Consider an arbitrary point $x \in M \setminus (K_{\frac{\delta}{2}}
\cup \{ x_1, ..., x_N \})$. \\
If $d_g(x, x_i) >2 r_i$ for all $i=1,...,N$, then
\begin{eqnarray*}
d_g\big(x, \cup _{i=1}^{N} B_{r_i}(x_i)\big) & = & d_g \left (x,
B_{r_i}(x_i) \right )
\hspace{.3in} \text { for some } i\\
& =& d_g(x, x_i) - r_i\\
& > & \frac{1}{2} d_g(x, x_i) \\
& \geq &\frac{1}{2} d_g(x, \{ x_1, ..., x_N \}). 
\end{eqnarray*}
So 
$$
u(x) \leq C \left (d_g \left (x, K_{\frac{\delta}{2}} \cup \{ x_1, ...,
x_N \} \right ) \right )^{-\frac{2}{p-1}}.
$$
If $x \in B_{2r_i}(x_i)$ for some $ i \in \{1,...,N\}$, then $2 r_i >
d_g(x, x_i)$, i.e., $$r_i > \frac{1}{2}d_g(x, x_i) \geq \frac{1}{2}
d_g \left (x, K_{\frac{\delta}{2}}
\cup \{ x_1, ..., x_N \} \right ) . $$ In the coordinate system $y$ centered
around $x_i$ as in Lemma \ref{lem:iso},   
$$ \Bigg \|
u(x_i)^{-1}u \left (\frac{y}{u(x_i)^{\frac{p-1}{2}}} \right ) -
\left (1+\frac{K(x_i)}{n(n-2)} |y|^2 \right )^{-\frac{n-2}{2}} \Bigg \| <
\epsilon,  $$ therefore
\begin{eqnarray*}
u(x) & \leq & (1+ \epsilon) u(x_i)\\
& = &  (1+ \epsilon)R^{ \frac{2}{p-1}} r_i^{-\frac{2}{p-1}}
\hspace{.3in} \text{ since } r_i=\frac{R}{u(x_i)^{\frac{p-1}{2}}}\\ 
& \leq & C \bigg (d_g \left ( x, K_{\frac{\delta}{2}} \cup \{ x_1, ...,
x_N \} \right ) \bigg )^{-\frac{2}{p-1}}.
\end{eqnarray*}

We conclude that
\begin{prop}
\label{prop:iso}
Given $\epsilon >0, R >> 0$, there exists $C=C(\epsilon, R)$ such that
if $u$ is a solution of equation $(\ref{eq:main})$ and
$$\max _{x \in M} \Big( \big(d_g(x, K_{\frac{\delta}{2}})\big
)^{\frac{2}{p-1}}u(x) \Big) > C,$$
then there exists $\{ x_1, ..., x_N \} \subset M \setminus
K_{\frac{\delta}{2}}$ with $N $ depending on $u$, and

\begin{itemize}

\item Each $x_i$ is a local maximum of $u$ and the geodesic balls $\{
B_{r_i} (x_i)\}$ are disjoint. 

\item $|\frac{n+2}{n-2}-p | < \epsilon$ and in the coordinate system $y$ so
chosen that $z=\frac{y}{u(x_i)^{\frac{p-1}{2}}}$ is the geodesic normal coordinate system centered at $x_i$, we have
 
$$ \Bigg \| u(x_i)^{-1}u\left (\frac
{y}{u(x_i)^{\frac{p-1}{2}}} \right ) -\bar{v}(y) \Bigg \|_{C^2(B_{2R}(0))} <
\epsilon$$ on the ball $B_{2R}(0) \subset \mathbf{R}^n(y)$, where
$$\bar{v} (y)=\left (1+
\frac{K(x_i)}{n(n-2)}|y|^2 \right )^{-\frac{n-2}{2}}.$$ 

\item There exists $C=C(\epsilon, R)$ such that $$u(x) \leq C \left ( d_g(x,
\overline{K_{\frac{\delta}{2}} \bigcup \{ x_1, ..., x_N \}}) \right )^{-\frac{2}{p-1}}.        $$

\end{itemize} 

\end{prop}

Before we proceed with the proof, we need to give two definitions.

\begin{defn}
\label{defn:bupt}
We call a point $\bar{x}$ on a manifold $M$ a {\bf blow-up point} of the sequence $\{ u_i \}$ if $\bar{x}=\lim _{i
\rightarrow \infty} x_i$ for some $\{ x_i \} \subset M$ and $u_i(x_i) \rightarrow \infty$.
\end{defn}

\begin{defn}
\label{defn:iso}
Let $\{ u_i \}$ be a sequence of functions satisfying $\Delta _{g_i} u_i + K_i
u_i^{p_i}=0$ on some manifold $M$ where the metrics $g_i$ converge to some metric
$g_0$. A point $\bar{x} \in M$ is called an {\bf isolated blow-up
point} of $\{ u_i \}$  if there exist local maximum points $x_i$ of
$u_i$ and a fixed radius $r_0>0$ such that 
\begin{enumerate}
\item $x_i \rightarrow \bar{x}$.
\item $u_i(x_i) \rightarrow \infty$.
\item $u_i(x) \leq C \left (d_{g_i} (x, x_i)\right )^{-\frac{2}{p_i-1}} \hspace{.1in} \forall x
\in B_{r_0}(x_i)$ for some constant $C$ independent of $i$.
\end{enumerate}

\end{defn}

Now we are going to prove that $u$ is uniformly bounded on
$M \setminus K_{\delta}$. Suppose it is not, then there are sequences 
$\{ u_i \}, \{ p_i \} $ and $\{ K_i \}$ such that $$ \Delta _{g_i} u_i + K_i
u_i^{p_i}=0 \hspace{.3in} \text{ and } \hspace{.2in} \max _{M \setminus K_{i,\delta}} u_i
\rightarrow \infty \,\, \text{ as } i \rightarrow \infty$$  where
$K_{i, \delta}:=\{x \in M: K_i(x) < \delta\}$.

By an argument similar to that in section \ref{section:smallK}, $d_g(\partial K_{i,\delta},
\partial K_{i,\frac{\delta}{2}})$ is bounded below by some constant depending only
on $M, g,n, \zeta, \Lambda, C_K$, so $\max _{M \setminus K_{i,\delta}} (
[d_g(x, K_{i,\frac{\delta}{2}})]^{\frac{2}{p_i-1}}u_i) \rightarrow \infty$
as $i \rightarrow \infty $. Thus for fixed $\epsilon >0$ and $R >> 0
$  we can apply Proposition \ref{prop:iso} to each $u_i$ and
find $x_{1,i},..., x_{N(i),i}$ such that 
\begin{equation}
\label{locmax}
\text{ each } x_{j,i} \,\,
(1 \leq j \leq N(i) ) \text{ is a local maximum of } u_i;
\end{equation}

\begin{equation}
\label{disjoint}
\text{the balls } B_{\frac{R}{u(x_{j, i})^{\frac{p_i -1}{2}}}}(x_{j, i}) \text{ are disjoint};   
\end{equation}

\noindent
and for coordinates $y$ centered at $ x_{j,i}$ such that
$\frac{y}{u(x_{j,i})^{\frac{p_i-1}{2}}}$ is the geodesic normal coordinates,
\begin{equation}
\label{sphe}
\Bigg \| u_i(x_{j,i})^{-1}u_i \left (\frac
{y}{u_i(x_{j,i})^{\frac{p_i-1}{2}}}\right ) - \left (1+
\frac{K_i(x_{j, i})}{n(n-2)}|y|^2 \right )^{-\frac{n-2}{2}} \Bigg  \|_{C^2(B_{2R}(0))} <\epsilon .
\end{equation}

Let $\sigma _i= \min \{ d_g(x_{\alpha, i}, x_{\beta ,i}): \alpha \neq
\beta, 1 \leq \alpha, \beta \leq N(i) \}$. Without lost of generality we
can assume $\sigma _i= d_g(x_{1, i}, x_{2 ,i})$. There are two
possibilities which could happen.

\noindent
{\bf Case I}: $\sigma _i \geq \sigma>0$. \\
Then the points $x_{j,i}$ have isolated
limiting points $x_1, x_2, ...$, which are isolated blow-up points of $\{
u_i \}$ as defined above.  

\noindent
{\bf Case II}: $\sigma _i \rightarrow 0$.\\ 
Then we rescale the coordinates to make
the minimal distance to be 1: let $y=\sigma _i ^{-1}z$ where $z$ is the
geodesic normal coordinate system centered at $x_{1,i}$. We also rescale the function by defining
$$v_i(y)=\sigma _i ^{\frac{2}{p_i-1}}u_i(\sigma _i y).$$
$v_i$ satisfies 
$$
\Delta _{g^{(i)}} v_i + K_i(\sigma _i y)v_i^{p_i}=0
$$
where the metric
$g^{(i)}(y)=g_{\alpha \beta}(\sigma _i y)dy^{\alpha}dy^{\beta}$.

Let
$y_{j,i}$ be the coordinate corresponding to $x_{j,i}$ then the
distance between any two points in $\{y_{1,i}, ... , y_{N(i),i}\}$ is at least $1$. Let $\{ y_1, y_2, ...  \}$ be the limiting points of $\{ y_{j,i}\}$.

Let $\Omega$ be any compact subset of $\mathbf{R}^n(y) \setminus \{
y_1, y_2, ... \}$. Because we have proved that $u_i$ is uniformly
bounded on $K_{i,\delta}$, we must have the maximum points $\{x_{j,i}\} \subset M
\setminus K_{i,\delta}$ and therefore $$d_g(x_{j,i},
K_{i,\frac{\delta}{2}}) \geq d_g(\partial K_{i,\delta}, \partial K_{i,\frac{\delta}{2}}) \geq
C.$$ This means for $y \in \Omega$ and $x=\exp _{x_{1,i}}(\sigma_i y) \in M$, when $i$ is
large enough, $$d_g\left (x, K_{i,\frac{\delta}{2}} \cup
\{x_{1,i},...,x_{N(i),i}\}\right)=d_g \left (x,
\{x_{1,i},...,x_{N(i),i}\} \right).$$
Then by Proposition \ref{prop:iso} and the fact that $g^{(i)}$ converges
to the Euclidean metric on $\Omega$ ,
\begin{eqnarray*}
v_i(y) & = & \sigma _i^{\frac{2}{p_i-1}}u_i(\sigma _i y)\\
& \leq &\sigma _i^{\frac{2}{p_i-1}} C(\epsilon, R)\left ( d_g \left (x,
\{x_{1,i},...,x_{N(i),i}\}\right )\right )^{-\frac{2}{p_i-1}} \\
&=&C(\epsilon, R)\left ( d_{g^{(i)}}\left (y, \{y_{1,i},...,y_{N(i),i}\}
\right )
\right )^{-\frac{2}{p_i-1}}\\
& \leq & C(\epsilon, R, \Omega) 
\end{eqnarray*}

\noindent
We also know that $K_i$ has uniform $C^3$ bound, so it
converges in $C^2$ norm to some function $K$. 
Then by the standard elliptic estimates $\{v_i\}$
converge on $\Omega$ to some function $v$ in $C^2$-norm  which
satisfies $\Delta v+K(x_1)v^{\frac{n+2}{n-2}}=0$ on $\Omega$, where
$x_1$ is the limit point of $\{x_{1,i}\}$ and $K(x_1) \geq \delta >0$. 
Since $\Omega$ is arbitrary, $\Delta v+K(x_1)v^{\frac{n+2}{n-2}}=0$ on
$\mathbf{R}^n(y) \setminus \{y_1,y_2,...\}$.


If none of the points $0=y_1, y_2,...$ is a blow-up point for
$\{v_i\}$, then $\{v_i \}$ is uniformly bounded near each of those
points. Thus the
convergence $v_i \rightarrow v$ holds near each of those points and
therefore  $\Delta v+K(x_1)v^{\frac{n+2}{n-2}}=0$ on $\mathbf{R}^n(y)$. By the definition of $v_i$ and the choice of $\sigma
_i$, $x_{1,i}$ and $x_{2,i}$,
\begin{eqnarray*}
v_i(0)&=&\sigma _i ^{\frac{2}{p_i-1}} u_i(x_{1,i})\\
& > & \left (\frac {R}{u_i(x_{1,i})^{\frac{p_i-1}{2}}} \right )^{\frac{2}{p_i-1}}
u_i(x_{1,i}) \hspace{.5in} \text{(by (\ref{disjoint}))}\\
&=& R^{\frac{2}{p_i-1}}.
\end{eqnarray*} 
So $v(0)$ is bounded below away from $0$ and hence $v
\neq 0$. By similar argument as before we then know $v$ can be
expressed as $$ v(y)=\left (\frac{\lambda
K(x_1)^{-\frac{1}{2}}}{1+\frac{1}{n(n-2)}\lambda ^2 |y-\bar{y}|^2
}\right )^{\frac{n-2}{2}} $$ for some $\lambda >0$ and $\bar{y} \in
\mathbf{R}^n(y)$. It implies that $v$ can only have
one critical point. But we know each $v_i$ has at least two critical
points $0$ and $y_{2,i}$, $|y_{2,i}|=1$, and $v_i$ is $C^2$-close to $v$, so $v$ must
also have at least two critical points. This is a contradiction. Thus $\{v_i\}$
has at least one blow-up point, without lost of generality we can
assume it to be $0$.

If there are other blow-up points besides $0$, they are at least
distance $1$ apart. For any $|y| \leq \frac{1}{2}$, the corresponding
$x=\exp_{x_{1,i}}(\sigma _i y), x_{j,i}=\exp_{x_{1,i}}(\sigma _i y_{j,i}) \in M$ satisfy  $d_g(x, \{x_{1,i}, ... , x_{N(i), i} \}) = d_g(x, x_{1,i})$. So by Proposition \ref{prop:iso}
\begin{eqnarray*} 
v_i(y) & = & \sigma _i ^{\frac{2}{p_i-1}} u_i(\sigma _i y)\\
 & = & \sigma _i
^{\frac{2}{p_i-1}} u_i(x)\\
& \leq & \sigma _i ^{\frac{2}{p_i-1}} C(\epsilon, R) d_g(x, x_{1,i})^{-\frac{2}{p_i-1}}\\
& = & C(\epsilon, R) d_{g^{(i)}}(y, 0)^{-\frac{2}{p_i-1}}
\end{eqnarray*}

Therefore we have reduced Case II to the following case:\\
There is a sequence of functions $\{v_i\}$, each satisfies  
\begin{equation}
\label{eq:rscase2}
\Delta _{g^{(i)}} v_i + K_i(\sigma _i y)v_i^{p_i}=0
\end{equation}
where $g^{(i)}(y)=g_{\alpha \beta}(\sigma _i y)dy^{\alpha}dy^{\beta}$ converges to the
Euclidean metric on compact subset of $\mathbf{R}^n(y)$. The sequence
$\{ v_i \}$ has isolated blow-up point(s) $\{0,...\}$. 

In the following sections we are going to show neither Case I nor Case
II can happen for $n=3,4$.

\section{Simple blow-up and Related Estimates}
\label{section:simple}

In this section we are going to diverge from the proof of the
compactness theorems temporarily. We will analyze the phenomenon of
simple blow-up and obtain some estimates which are important in the rest of the proof.

\begin{defn}
\label{defn:simple}
$x_0$ is called a {\bf simple blow-up point} of $\{u_i\}$ if it is an isolated blow
up point and there exists $r_0 >0$ independent of $i$ such that
$\bar{w}_i(r)$ has only one critical point for $r \in (0, r_0)$, where $\bar{w}_i(r)=\Vol (S_r)^{-1} \int _{S_r} |z|^{\frac{2}{p_i-1}}u_i(z)d \Sigma _g$, and $z$ is the local coordinate system centered at each $x_i$.
\end{defn} 

We are going to derive some estimates of $u_i$ near a simple
blow-up point.

The first lemma actually only requires $x_0$ to be an isolated blow-up
point.

\begin{lemma}
\label{lemma:simpleharnack}
If $x_0=\lim _{i \rightarrow \infty} x_i$ is an isolated blow-up point
of $\{u_i\}$ which satisfies $\Delta _{g^{(i)}}u_i+K_iu_i^{p_i}=0$, then there exists a constant $C$ independent of $i$ and
$r$ such that $\max _{\partial B_r(x_i)} u_i(x) \leq C \min _{\partial
B_r(x_i)} u_i(x)$
for any $0<r \leq r_0$ where $r_0$ is the fixed radius as in
Definition {\em \ref{defn:iso}}.
\end{lemma}

\pf
Let $z$ be the coordinates centered at each $x_i$ and $y=\frac{z}{r}$.
Define $\hat{u}_i(y)=r^{\frac{2}{p_i-1}}u_i(ry)$. Let
$\hat{g}_i(y)=g^{(i)}_{\alpha \beta}(ry)dy^{\alpha}dy^{\beta}$. Then
$$
\Delta _{\hat{g}(y)}\hat{u}_i+K_i(ry)\hat{u}_i^{p_i}=0.$$
Since $
u_i(z) \leq C |z|^{-\frac{2}{p_i-1}}  $ for $|z| \leq r_0$, we know
$|y|^{\frac{2}{p_i-1}}\hat{u}_i(y) \leq C$ for $|y| \leq \frac{r_0}{r}$.
In particular, when $|y|=1$, $\hat{u}_i(y) \leq C$. Then we can apply
the 
standard Harnack inequality to get $\max _{|y|=1} \hat{u}_i(y) \leq C
\min _{|y|=1} \hat{u}_i(y)$ for some constant $C$ independent of $i$ and
$r$, so $$\max _{|z|=r} u_i(z) \leq C
\min _{|z|=r} u_i(z).$$
\stop

\begin{prop}
\label{prop:simpleestimates}
Let $x_0=\lim _{i \rightarrow \infty}x_i$ be a simple blow-up point of
$u_i$ with $p_i \rightarrow \frac{n+2}{n-2} $. Let $z$ be the geodesic coordinates centered at each
$x_i$. There exist constants $\bar{r} \leq r_0$ and $C$ independent of
$i$ such that
\begin{itemize}
\item if \,\, $0
\leq |z| \leq \bar{r}$, \hspace{.1in} then $$u_i(z) \geq C u_i(x_i)\left (1+ \frac{K_i(x_i)}{n(n-2)}
u_i(x_i)^{\frac{4}{n-2}}|z|^2 \right )^{-\frac{n-2}{2}}$$ 
\item  if \,\,  $0
\leq |z| \leq \frac{R}{u_i(x_i)^{\frac{p_i-1}{2}}}$, \hspace{.1in}
then $$u_i(z) \leq C u_i(x_i) \left (1+ \frac{K_i(x_i)}{n(n-2)}
u_i(x_i)^{p_i-1}|z|^2 \right )^{-\frac{n-2}{2}}$$ 
\item if \,\, $\frac{R}{u_i(x_i)^{\frac{p_i-1}{2}}} \leq |z| \leq \bar{r}$,
\hspace{.1in} then \hspace{.1in} $u_i(z) \leq Cu_i(x_i)^{t _i}|z|^{-l_i}$\\
 where $l_i$, $t_i$ are chosen such that 
$l_i \rightarrow \frac{6(n-2)}{7}$ when $i \rightarrow \infty$,
 and $t _i =1-\frac{(p_i-1)l_i}{2}$.
\end{itemize} 
\end{prop}

\pf
By Proposition \ref{prop:iso}, when $0 \leq |z| \leq \frac{R}{u_i(x_i)^{\frac{p_i-1}{2}}}$,
\begin{eqnarray*}
& & (1-\epsilon) \frac{u_i(x_i)}{\left (1+ \frac{K_i(x_i)}{n(n-2)}
u_i(x_i)^{p_i-1}|z|^2 \right )^{\frac{n-2}{2}}} \\
&  \leq & u_i(z)\\
& \leq & (1+\epsilon) \frac{u_i(x_i)}{\left (1+ \frac{K_i(x_i)}{n(n-2)}
u_i(x_i)^{p_i-1}|z|^2 \right )^{\frac{n-2}{2}}}.  
\end{eqnarray*}

\noindent
Since $u_i(x_i) >1$ and $p_i-1 \leq \frac {4}{n-2}$,
\begin{eqnarray*}
& &  u_i(x_i)\left (1+ \frac{K_i(x_i)}{n(n-2)}
u_i(x_i)^{p_i-1}|z|^2 \right )^{-\frac{n-2}{2}} \\
& \geq & u_i(x_i)\left (1+ \frac{K_i(x_i)}{n(n-2)}
u_i(x_i)^{\frac{4}{n-2}}|z|^2\right )^{-\frac{n-2}{2}}.
\end{eqnarray*}
So we only need to find the upper and lower bounds for $u_i(z)$ when $\frac{R}{u_i(x_i)^{\frac{p_i-1}{2}}} \leq |z| \leq \bar{r}$.

First the lower bound. 

Let $G_i$ be the Green's function of $\Delta _{g_i}$ which is singular at
$0$ and $G_i =0$ on $\partial B_{\bar{r}}$. Then $G_i(z)=|z|^{2-n}+R_i(z)$
where $\lim _{|z| \rightarrow 0} |z|^{n-2}R_i(z)=0$. Since $g_i$
converges uniformly to $g_0$, 
there exist constants $C_1$ and $C_2$ independent of $i$ such that
$$C_1 \leq 1+|z|^{n-2}R_i(z)=|z|^{n-2}G_i(z)\leq C_2 \hspace{.1in} \text{ for
} |z| \leq \bar{r},$$ i.e. $C_1 |z|^{2-n} \leq G_i(z) \leq C_2 |z|^{2-n}$. 

When $|z|=Ru_i(x_i)^{-\frac{p_i-1}{2}}$,
\begin{eqnarray*}
u_i(z) & \geq & (1-\epsilon) \frac{u_i(x_i)}{\left (1+ \frac{K_i(x_i)}{n(n-2)}
u_i(x_i)^{p_i-1}|z|^2 \right )^{\frac{n-2}{2}}}\\
& = & (1-\epsilon)\frac{u_i(x_i)}{\left (1+ \frac{K_i(x_i)}{n(n-2)}
R^2 \right )^{\frac{n-2}{2}}}
\end{eqnarray*}
and 
\begin{eqnarray*}
u_i(x_i)^{-1} G_i(z) & \leq &  C_2 |z|^{2-n}u_i(x_i)^{-1}\\
 & = & C_2R^{2-n}u_i(x_i)^{\frac{(n-2)(p_i-1)}{2}-1}\\
& \leq & C_2 R^{2-n}u_i(x_i) \hspace{.3in} (\text{ since }
\frac{(n-2)(p_i-1)}{2}-1 \leq 1 \,\,).
\end{eqnarray*}
For $R>>0$,
\begin{eqnarray*}
R^{n-2}\left (1+ \frac{K_i(x_i)}{n(n-2)}
R^2\right )^{-\frac{n-2}{2}} & = & \left (R^{-2}+
\frac{K_i(x_i)}{n(n-2)} \right)^{-\frac{n-2}{2}}\\
& \geq & C(n, C_K).
\end{eqnarray*}
Therefore $u_i(z) \geq C u_i(x_i)^{-1}G_i(z)$ for some constant $C$
independent of $i$ when $|z|=Ru_i(x_i)^{-\frac{p_i-1}{2}}$.\\
For that constant $C$, $u_i(z) \geq C u_i(x_i)^{-1}G_i(z)=0$ when $|z|=\bar{r}$.\\
Then since $$\Delta \left (u_i(z)- C  u_i(x_i)^{-1}G_i(z)\right )=\Delta
u_i(z)=-K_iu_i(z)^{p_i}<0$$ on $B_{\bar{r}}$, by the maximal principle
$$u_i(z)-Cu_i(x_i)^{-1}G_i(z)>0, \hspace{.3in} \text{ i.e., } \hspace{.3in} u_i(z) > Cu_i(x_i)^{-1}G_i(z)
$$ when $\frac{R}{u_i(x_i)^{\frac{p_i-1}{2}}} \leq |z| \leq \bar{r}$.

\noindent
Because $u_i(x_i)^{-1}G_i(z) \geq C_1|z|^{2-n}u_i(x_i)^{-1}$, we now need
to compare $|z|^{2-n}u_i(x_i)^{-1}$ with $u_i(x_i)\cdot\left (1+ \frac{K_i(x_i)}{n(n-2)}
u_i(x_i)^{\frac{4}{n-2}}|z|^2 \right )^{-\frac{n-2}{2}}$ in order to get the
desired lower bound.

\begin{eqnarray*}
 & & u_i(x_i)^2|z|^{n-2} \left (1+ \frac{K_i(x_i)}{n(n-2)}
u_i(x_i)^{\frac{4}{n-2}}|z|^2 \right )^{-\frac{n-2}{2}} \\
& \leq &  u_i(x_i)^2 \left (\frac{K_i(x_i)}{n(n-2)}
u_i(x_i)^{\frac{4}{n-2}} \right )^{-\frac{n-2}{2}}\\
& \leq & C  
\end{eqnarray*} 
where the constant $C$ is independent of $i$ because $K_i(x_i) \geq
\delta$ which doesn't depend on $i$. Therefore $$|z|^{2-n}u_i(x_i)^{-1}
\geq C u_i(x_i) \left (1+ \frac{K_i(x_i)}{n(n-2)}
u_i(x_i)^{\frac{4}{n-2}}|z|^2 \right )^{-\frac{n-2}{2}}, $$
which then implies that 
$$u_i(z) \geq C u_i(x_i) \left (1+ \frac{K_i(x_i)}{n(n-2)}
u_i(x_i)^{\frac{4}{n-2}}|z|^2 \right )^{-\frac{n-2}{2}}$$ 
when $\frac{R}{u_i(x_i)^{\frac{p_i-1}{2}}} \leq |z| \leq \bar{r}$.

Next the upper bound. We are going to apply the same strategy of
constructing a comparison function and using the maximal principle.

Define $\mathcal{L}_i\varphi := \Delta _{g_i} \varphi + K_i
u_i(z)^{p_i-1}\varphi$. By definition $\mathcal{L}_iu_i=0$. Let $M_i=\max
_{\partial B_{\bar{r}}} u_i$ and $m_i=\min
_{\partial B_{\bar{r}}} u_i$. Let  $C_i=(1+\epsilon)
\left (\frac{K_i(x_i)}{n(n-2)} \right )^{-\frac{n-2}{2}}$. $C_i$ is bounded above
and below by constants only depending on $\epsilon, n, C_K$ and $\delta$. Consider the function
$$M_i(\bar{r}^{-1}|z|)^{-n+2+l_i}+C_iu_i(x_i)^{t _i}|z|^{-l_i}.$$

When $|z|=\frac{R}{u_i(x_i)^{\frac{p_i-1}{2}}}$,
 
\begin{eqnarray}
u_i(z) & \leq & (1+\epsilon) \frac{u_i(x_i)}{\left (1+ \frac{K_i(x_i)}{n(n-2)}
u_i(x_i)^{p_i-1}|z|^2 \right )^{\frac{n-2}{2}}} \nonumber \\
&=& (1+\epsilon) \frac{u_i(x_i)}{\left (1+ \frac{K_i(x_i)}{n(n-2)}R^2
\right )^{\frac{n-2}{2}}} \nonumber \\
& \leq & C_iu_i(x_i)R^{-l_i} \hspace{.45in} \text{ (because } l_i < n-2) \nonumber \\
& = & C_i u_i(x_i)^{t _i}|z|^{-l _i} \hspace{.35in}\text{
(by the choice of } t _i). \nonumber
\end{eqnarray}
\noindent
When $|z|=\bar{r}$, by the definition of $M_i$, $u_i(z) \leq M_i=M_i(\bar{r}^{-1}|z|)^{-n+2+l_i}.$ \\
So on $\{|z|=\bar{r} \} \cup \{|z|=Ru_i(x_i)^{-\frac{p_i-1}{2}}
\}$,$$u_i(z) \leq M_i(\bar{r}^{-1}|z|)^{-n+2+l_i}+C_iu_i(x_i)^{t _i}|z|^{-l_i}.$$

In the Euclidean coordinates, $\Delta |z|^{-l_i}=-l_i(n-2-l_i)|z|^{-l_i-2}
$ and $\Delta |z|^{-n+2+l_i}=-l_i(n-2-l_i)|z|^{-n+l_i}
$. Since $z$ is the geodesic normal coordinates, when $\bar{r}$ is sufficiently small, $g_0$
and $g_i$ are close to the Euclidean metric. Then
when $i$ is large enough 
\begin{equation}
\Delta _{g_i} |z|^{-l_i} \leq -\frac{1}{2}l_i(n-2-l_i)|z|^{-l_i-2}
\end{equation}
and 
\begin{equation}
\label{eq:|z|^-n+2+l_i}
\Delta _{g_i} |z|^{-n+2+l_i} \leq
-\frac{1}{2}l_i(n-2-l_i)|z|^{-n+l_i}. 
\end{equation}
Thus
\begin{eqnarray*}
& & \mathcal{L}_i(C_iu_i(x_i)^{t _i}|z|^{-l_i}) \\
&=& C_i u_i(x_i)^{t
_i}\Delta _{g_i}|z|^{-l_i}+C_iu_i(x_i)^{t _i}K_iu_i(z)^{p_i-1}|z|^{-l_i} \\
& \leq & -Cl_i(n-2-l_i)u_i(x_i)^{t _i}|z|^{-l_i-2}+C'u_i(x_i)^{t _i}u_i(z)^{p_i-1}|z|^{-l_i}
\end{eqnarray*}
for some constants $C, C'$ independent of $i$.

\noindent
The upper bound on $u_i(z)$ when $|z| =
Ru_i(x_i)^{-\frac{p_i-1}{2}}$ and Lemma \ref{lemma:simpleharnack} implies
that
\begin{eqnarray*}
\bar{u}_i \left (Ru_i(x_i)^{-\frac{p_i-1}{2}} \right ) & \leq & 
\frac{(1+\epsilon)u_i(x_i)}{ \left [1+ \frac{K_i(x_i)}{n(n-2)}
u_i(x_i)^{p_i-1} \left ( Ru_i(x_i)^{-\frac{p_i-1}{2}} \right )^2 \right ]^{\frac{n-2}{2}}}\\
& \leq & Cu_i(x_i)R^{2-n}.
\end{eqnarray*}
Since $x_0$ is a simple point of blow-up, $r^{\frac{2}{p_i-1}}\bar{u}
_i(r)$ is decreasing from $Ru_i(x_i)^{-\frac{p_i-1}{2}}$ to $\bar{r}$,
which implies
\begin{eqnarray*}
|z|^{\frac{2}{p_i-1}} \bar{u}_i(|z|) & \leq &
\left (Ru_i(x_i)^{-\frac{p_i-1}{2}} \right )^{\frac{2}{p_i-1}} \cdot \bar{u}_i
\left (Ru_i(x_i)^{-\frac{p_i-1}{2}} \right )\\
& \leq & CR^{\frac{2}{p_i-1}+2-n}.
\end{eqnarray*}
Thus by Lemma \ref{lemma:simpleharnack} again  
\begin{equation}
\label{eq:u_i^p_i-1|z|^-l_i}
u_i(z)^{p_i-1} \leq
C|z|^{-2}R^{2-(n-2)(p_i-1)}
\end{equation}
and hence 
$$
u_i(z)^{p_i-1} |z|^{-l_i}\leq
C|z|^{-2-l_i}R^{2-(n-2)(p_i-1)}.
$$
So we know
\begin{eqnarray*}
& & \mathcal{L}_i \left (C_iu_i(x_i)^{t _i}|z|^{-l_i} \right )\\
 &   \leq & \left
(-Cl_i(n-2-l_i)+C' R^{2-(n-2)(p_i-1)} \right )  u_i(x_i)^{t _i}|z|^{-l_i-2} 
\end{eqnarray*}
By our choice of $l_i$, \,\, $l_i(n-2-l_i)$ is always
bounded below by some positive constant independent of $i$. When $i$
is sufficiently large, $2-(n-2)(p_i-1) < 0$, we can choose $R$ big enough
such that $-Cl_i(n-2-l_i)+C' R^{2-(n-2)(p_i-1)}<0 $, hence $\mathcal{L}_i(C_iu_i(x_i)^{t _i}|z|^{-l_i}) <0$.

\noindent
Similarly,

\begin{eqnarray*}
& & \mathcal{L}_i \left (M_i(\bar{r}^{-1}|z|)^{-n+2+l_i} \right )\\
 & = & M_i \bar{r}^{n-2-l_i}
\Delta _{g_i} |z|^{-n +2 + l_i} +M_i  \bar{r}^{n-2-l_i} K_i
u_i(z)^{p_i-1}|z|^{-n+2+l_i}\\
& \leq & 
-\frac{1}{2}l_i(n-2-l_i)M_i \bar{r} ^{n-2-l_i}|z|^{-n+l_i} \\
& & 
+K_i M_i \bar{r}^{n-2-l_i}R^{2-(n-2)(p_i-1)}|z|^{-n+l_i} 
\end{eqnarray*}
by equations (\ref{eq:|z|^-n+2+l_i}) and (\ref{eq:u_i^p_i-1|z|^-l_i}). We can choose $R$ large enough such that
$-\frac{1}{2}l_i(n-2-l_i) +K_iR^{2-(n-2)(p_i-1)} <0$ and hence
$$\mathcal{L}_i(M_i(\bar{r}^{-1}|z|)^{-n+2+l_i})<0. $$

Therefore when $Ru_i(x_i)^{-\frac{p_i-1}{2}} \leq |z|
\leq \bar{r}$,
$$\mathcal{L}_i \left (M_i(\bar{r}^{-1}|z|)^{-n+2+l_i}+C_iu_i(x_i)^{t
_i}|z|^{-l_i} \right ) <0. $$  Then by the maximal principle  $$u_i(z) \leq M_i(\bar{r}^{-1}|z|)^{-n+2+l_i}+C_iu_i(x_i)^{t
_i}|z|^{-l_i}.$$ 

By Lemma \ref{lemma:simpleharnack} and because $x_0$ is a simple blow-up point, for
$\frac{R}{u_i(x_i)^{\frac{p_i-1}{2}}} $ $\leq  \theta \leq \bar{r}$, 
\begin{eqnarray*}
\bar{r}^{\frac{2}{p_i-1}}M_i & \leq & \theta^{\frac{2}{p_i-1}}
\bar{u}_i(\theta) \\
& \leq &
\theta ^{\frac{2}{p_i-1}} \left (M_i(\bar{r}^{-1} \theta)^{-n+2+l_i}+C_iu_i(x_i)^{t
_i} \theta ^{-l_i} \right )\\
 & = & \bar{r}^{n-2-l_i} \theta ^{\frac{2}{p_i-1}-n+2+l_i}M_i + \theta
^{\frac{2}{p_i-1}} \cdot C_iu_i(x_i)^{t
_i} \theta ^{-l_i} 
\end{eqnarray*}
for some constant $C$ independent of $i$.

\noindent
When $i \rightarrow \infty$, $\frac{2}{p_i-1}-n+2+l_i \rightarrow \frac{5}{14}(n-2)>0$.

\noindent
Since $\frac{R}{u_i(x_i)^{\frac{p_i-1}{2}}} \to 0$, we can choose $\theta$ small enough (fixed, independent of $i$) to
absorb the first term on the right hand side of the above inequality
into the left hand side to get $M_i \leq 2 C_i \theta^{\frac{2}{p_i-1}-l_i}u_x(x_i)^{t_i}
\leq C u_i(x_i)^{t _i}.$  Therefore
\begin{eqnarray*}
u_i(z) & \leq & M_i(\bar{r}^{-1}|z|)^{-n+2+l_i}+C_iu_i(x_i)^{t
_i}|z|^{-l_i} \\
& \leq & M_i(\bar{r}^{-1}|z|)^{-l_i}+C_iu_i(x_i)^{t
_i}|z|^{-l_i} \\
& \leq & Cu_i(x_i)^{t
_i}|z|^{-l_i}
\end{eqnarray*}

\stop

\begin{prop}
\label{prop:u^delta}
If $x_0=\lim _{i \rightarrow \infty} x_i$ is a simple blow-up point
and $p_i \rightarrow \frac{n+2}{n-2}$.
Let $\delta _i = \frac {n+2}{n-2} - p_i $, then $\lim _{i \rightarrow
\infty} u_i(x_i)^{\delta _i} =1$.
\end{prop}

For the proofs of this proposition and theorems \ref{thm:main} and \ref{thm:cor}, we need
to use the following {\bf Pohozaev identity} as proved in \cite{S5}.

\begin{prop}
\label{prop:pohozaev}
{\bf (Schoen, 1988)} Let $(N,g)$ be an $n$-dimensional compact Riemannian
manifold with smooth boundary $\partial N$. Let $R$ denote the scalar
curvature function of $N$, and suppose $X$ is a conformal Killing
vector field on $N$. We then have the identity
\begin{equation}
\label{eq:pohozaev}
\int _N ( \mathcal{L}_X R) dv = \frac {2n}{n-2} \int _{\partial N}
(\Ric-n^{-1}Rg)(X, \nu) d \sigma,
\end{equation} 
where $\Ric(\cdot , \cdot) $ denotes the Ricci tensor of $N$ thought of
as a quadratic form on tangent vectors, $\mathcal{L} _X$ denotes the
Lie derivative, $\nu$ denotes the outward unit normal vector to
$\partial N$, $dv$ and $d \sigma$ are volume and surface measure (with
respect to $g$), respectively.
\end{prop}

We now prove Proposition \ref{prop:u^delta}.

\pf 
Choose the conformal coordinate $z$ centered at $x_i$ such that on the
small ball $|z| \leq \sigma$, $g$ can be written as $\lambda
(z)^{\frac{4}{n-2}} g_0$ where $g_0$ is the Euclidean metric. Choose the conformal Killing field
$X=\sum _{j=1}^{n} z^j \frac{\partial}{\partial z^j}$, we can apply
the Pohozaev identity to get
\begin{equation}
\label{eq:pohozaev-X(R)} 
\frac{n-2}{2n} \int _{B_{\sigma}} X(R_i) dv_{g_i} = \int _{\partial
B_{\sigma}} T_i(X, \nu _i) d \Sigma _i
\end{equation}

\noindent
where the notations are
\begin{eqnarray*}
g_i &  = &  u_i^{\frac{4}{n-2}}g \,\,  = \,\, (\lambda u_i)^{\frac{4}{n-2}}g_0, \\
R_i &  = & R(g_i) \,\, = \,\, c(n)^{-1} K_iu_i^{- \delta _i},\\
 dv_{g_i} & = & u_i^{\frac{2n}{n-2}} dv_g \,\, = \,\,  (\lambda u_i)^{\frac{2n}{n-2}}dz  ,\\
\nu _i & = & (\lambda u_i)^{-\frac{2}{n-2}}\sigma^{-1}\sum _{j} z^j
\frac{\partial}{\partial z^j} \\
& &  \text { is the unit outer
normal vector on } \partial B_{\sigma} \text{ with respect to } g_i,\\
d\Sigma _i & = & (\lambda u_i)^{\frac{2(n-1)}{n-2}} d \Sigma _{\sigma}
\\
& & 
\text{ where } d \Sigma _{\sigma} \text { is the surface element of the
standard } S^{n-1}(\sigma),\\
T_i & = & (n-2) (\lambda u_i)^{\frac{2}{n-2}} \left(\Hess \left ( (\lambda
u_i)^{-\frac{2}{n-2}} \right )-\frac{1}{n}\Delta \left ((\lambda
u_i)^{-\frac{2}{n-2}} \right ) g_0 \right )
\end{eqnarray*}
where $\Hess$ and $\Delta$ are taken with respect to the Euclidean metric $g_0.$

We are going to study the decay of both sides of (\ref{eq:pohozaev-X(R)}). 
 
Up to a constant the left hand side is
\begin{eqnarray}
& & c(n)\int _{B_{\sigma}} X(R_i) dv_{g_i} \nonumber \\
& = &\int _{B_{\sigma}} X(K_i
u_i^{-\delta _i}) (\lambda u_i)^{\frac{2n}{n-2}} dz \nonumber\\
& = & \int _{B_{\sigma}} X(K_i)u_i^{p_i+1} \lambda ^{\frac{2n}{n-2}} dz
- \delta _i \int  _{B_{\sigma}} K_i u_i^{p_i}X(u_i)
\lambda^{\frac{2n}{n-2}} dz \nonumber \\
& = & \int  _{B_{\sigma}}  |z| \frac {\partial K_i}{\partial r} u_i^{p_i+1}
\lambda ^{\frac{2n}{n-2}}dz 
 +
\frac{\delta _i}{p_i+1} \int _{B_{\sigma}} r \frac{ \partial
K_i}{\partial r} \lambda ^{\frac{2n}{n-2}} u_i^{p_i+1} dz \nonumber \\
& & + \frac{\delta _i}{p_i+1} \int
_{B_{\sigma}} K_i u_i^{p_i+1} r \frac{ \partial \lambda
^{\frac{2n}{n-2}}}{\partial r} dz \nonumber \\
& & - \frac{\delta _i}{p_i+1} \bigg [ - \int
_{ B_{\sigma}} K_i u_i^{p_i+1}\lambda^{\frac{2n}{n-2}} \diver X \,\, dz \nonumber \\
& &     +    \int
_{\partial B_{\sigma}} K_iu_i^{p_i+1}\lambda^{\frac{2n}{n-2}}  X \cdot
\left (\frac{\sum
z^j \frac{\partial}{\partial z^j}}{\sigma} \right )d\Sigma _{\sigma}
\bigg ] \nonumber 
\end{eqnarray}
which can be further written as

\begin{eqnarray}
\label{eq:X(R)in-u^delta}
& = & \left(1+ \frac{\delta _i}{p_i+1}\right) \int  _{B_{\sigma}}  |z| \frac
{\partial K_i}{\partial r} u_i^{p_i+1}\lambda ^{\frac{2n}{n-2}}dz \nonumber \\
& &  + \frac{\delta _i}{p_i+1} \int
_{B_{\sigma}}|z| K_iu_i^{p_i+1}  \frac{ \partial \lambda
^{\frac{2n}{n-2}}}{\partial r} dz  \nonumber \\
& &  +\frac{\delta _i}{p_i+1}n  \int
_{ B_{\sigma}} K_iu_i^{p_i+1}\lambda^{\frac{2n}{n-2}} dz - \frac{\delta
_i}{p_i+1} \int
_{\partial B_{\sigma}} \sigma
K_iu_i^{p_i+1}\lambda^{\frac{2n}{n-2}}d\Sigma _{\sigma} .
\end{eqnarray}

By Proposition \ref{prop:simpleestimates}
\begin{eqnarray*}
 \int _{|z| \leq \frac{R}{u_i(x_i)^{\frac{p_i-1}{2}}} } |z|
u_i(z)^{p_i+1} dz & \leq & C
u_i(x_i)^{p_i+1}  \int _{|z| \leq \frac{R}{u_i(x_i)^{\frac{p_i-1}{2}}}
} |z| dz   \\
& \leq & C u_i(x_i)^{p_i+1-\frac{(n+1)(p_i-1)}{2}}\\
& = & C u_i(x_i) ^{-\frac{2}{n-2} + \frac{n-1}{2}\delta _i}.
\end{eqnarray*}
Also since $\lim _{i \rightarrow \infty} \big ( n-l_i(p_i+1)+1
\big ) = -\frac{5}{7}n+1 <0$,
\begin{eqnarray*}
& & \int _{\frac{R}{u_i(x_i)^{\frac{p_i-1}{2}}} \leq |z| \leq \sigma } |z|
u_i(z)^{p_i+1} dz \\
 & \leq & C \int
_{\frac{R}{u_i(x_i)^{\frac{p_i-1}{2}}} \leq |z| \leq \sigma } |z|
\left ( u_i(x_i)^{t _i}|z|^{-l_i} \right )^{p_i+1} \\
            & \leq & C u_i(x_i)^{ t _i(p_i+1) -\frac{p_i-1}{2}\left
(n-l_i(p_i+1)+1\right )} \\
& = & C u_i(x_i)^{p_i+1-\frac{(n+1)(p_i-1)}{2}}\\
& & (\text{by the definition of } l_i \text{ and } t _i)\\
 & = & C u_i(x_i) ^{-\frac{2}{n-2} +\frac{n-1}{2}\delta_i}.
\end{eqnarray*}
So 
\begin{equation}
\label{eq:int-|z|u_i^p_i+1}
\int _{|z| \leq \sigma } |z| u_i(z)^{p_i+1} dz \leq C u_i(x_i)
^{-\frac{2}{n-2} + \frac{n-1}{2}\delta _i}
\end{equation}
and hence the first term in (\ref{eq:X(R)in-u^delta}) decays in the order of $u_i(x_i)
^{-\frac{2}{n-2} + \frac{n-1}{2}\delta_i}$ and the second term decays
even faster than that since $\delta _i \rightarrow 0$.\\
By Proposition \ref{prop:simpleestimates}, on $\partial B_{\sigma}$, $u_i$ decays in the
order of $u_i(x_i)^{t _i}$, so the fourth term in (\ref{eq:X(R)in-u^delta}) decays at
least in
the order of $u_i(x_i)^{t _i(p_i+1)}$.\\ The third term
$$\frac{\delta _i}{p_i+1}n  \int
_{ B_{\sigma}} K_iu_i^{p_i+1}\lambda^{\frac{2n}{n-2}} dz \geq C \delta
_i   \int _{ B_{\sigma}} u_i^{p_i+1} dz .$$
When $|z| \leq \frac{R}{u_i(x_i)^{\frac{p_i-1}{2}}}$, 
\begin{eqnarray*}
 u_i(z)& \geq & (1-\epsilon) \frac{u_i(x_i)}{\left (1+ \frac{K_i(x_i)}{n(n-2)}
u_i(x_
i)^{p_i-1}|z|^2 \right )^{\frac{n-2}{2}}}\\
& \geq & (1-\epsilon)
\frac{u_i(x_i)}{\left (1+ \frac{K_i(x_i)}{n(n-2)}R^2 \right )^{\frac{n-2}{2}}}\\
& \geq &  Cu_i(x_i),
\end{eqnarray*}
thus 
\begin{eqnarray}
\label{eq:lower-int-u^p+1}
 \int _{ B_{\sigma}} u_i^{p_i+1} dz & \geq & \int _{|z| \leq
\frac{R}{u_i(x_i)^{\frac{p_i-1}{2}} } } u_i^{p_i+1} dz \nonumber \\
& \geq & C u_i(x_i)^{p_i+1 - \frac{n}{2}(p_i-1)} \nonumber \\
& = & C u_i(x_i)^{ \frac{n-2}{2} \delta _i} \nonumber \\
& \geq & C.
\end{eqnarray}
So the third term is bounded below by $C \delta _i$.

Next we are going to study the decay of the right hand side of (\ref{eq:pohozaev-X(R)}).
\begin{eqnarray}
\label{eq:u^delta-Tterm}
& & \int _{\partial
B_{\sigma}} T_i(X, \nu _i) d \Sigma _i \nonumber \\
 & = & \int _{\partial
B_{\sigma}} (n-2) (\lambda u_i)^{\frac{2}{n-2}}
\bigg [ \Hess 
  \left ((\lambda u_i)^{-\frac{2}{n-2}} \right ) 
  \Big (r \frac {\partial}{\partial r}, (\lambda
  u_i)^{-\frac{2}{n-2}}\sigma^{-1} r \frac {\partial}{\partial r}
  \Big )  \nonumber \\  
& &    -\frac{1}{n}\Delta 
  \left ((\lambda u_i)^{-\frac{2}{n-2}} \right) 
  \left < r \frac {\partial}{\partial r},  (\lambda
  u_i)^{-\frac{2}{n-2}}\sigma^{-1} r \frac {\partial}{\partial r} \right > 
\bigg ] 
(\lambda u_i)^{\frac{2(n-1)}{n-2}}  d \Sigma _{\sigma} \nonumber \\
& & (\text{where } <\cdot, \cdot> \text { is the Euclidean metric})
\nonumber \\
& = & (n-2) \int _{\partial B_{\sigma}} \bigg [\sigma ^{-1}\Hess 
\left ( (\lambda u_i)^{-\frac{2}{n-2}} \right ) 
\left  (r \frac {\partial}{\partial r}, r \frac {\partial}{\partial r}
\right ) \nonumber \\
& & -\frac{\sigma}{n}\Delta \left ((\lambda
u_i)^{-\frac{2}{n-2}} \right ) \bigg ] (\lambda
u_i)^{\frac{2(n-1)}{n-2}} d \Sigma _{\sigma} \nonumber \\
& = & (n-2)
\int _{\partial B_{\sigma}} 
 \sigma ^{-1}  \bigg [-\frac {2 }{n-2} (\lambda u_i) \sum _{j,k}z^j z^k 
\frac{\partial}{\partial z^k} \frac{\partial}{\partial z^j} (\lambda
u_i)  \nonumber \\
& & +  \frac{2n}{(n-2)^2} \sum_{j,k} z^j z^k
\frac{\partial (\lambda u_i)}{\partial z^k} \frac{\partial( \lambda
u_i)}{\partial z^j} \bigg ]  - \sigma \cdot \\
& &\bigg [ -\frac{2 }{n(n-2)} (\lambda u_i) \sum _{j} \frac{\partial ^2 (\lambda u_i)}{(\partial
z^j)^2}  + \frac{2}{(n-2)^2} \sum _{j} \left(\frac{\partial (\lambda
u_i)}{\partial z^j} \right )^2 \bigg ]
d \Sigma _{\sigma} \nonumber
\end{eqnarray}

\noindent
On $\partial B _{\sigma} $, by Proposition \ref{prop:simpleestimates}, $u_i \leq C
u_i(x_i)^{t _i}$, so by the elliptic regularity theory \cite{GT} $\| u_i \| _{C^2(\partial B _{\sigma}) } \leq C
u_i(x_i)^{t _i}$. Thus we know (\ref{eq:u^delta-Tterm}) decays in the order of $u_i(x_i)^{2t _i}$.

Then by comparing the decay rate of both sides of (\ref{eq:pohozaev-X(R)}) 
\begin{equation}
\label{eq:upper-on-delta}
 \delta _i \leq C\left (u_i(x_i)^{-\frac{2}{n-2}+ \frac{n-1}{2}\delta
_i} +  u_i(x_i)^{2t _i}\right ).
\end{equation}
Thus 
$$
 \delta _i \ln u_i(x_i)\leq
C\left (u_i(x_i)^{-\frac{2}{n-2}+\frac{n-1}{2}\delta _i} +
u_i(x_i)^{2t _i} \right )\ln u_i(x_i).
$$

\noindent
By our choice of $l_i,$ $$  t _i = 1 -\frac{(p_i-1)l_i}{2} \rightarrow
-\frac{5}{7} <0$$
Since $u_i(x_i) \rightarrow \infty$, $\left
(u_i(x_i)^{-\frac{2}{n-2}+\frac{n-1}{2}\delta_i} + 
u_i(x_i)^{2t _i} \right)\ln u_i(x_i) \rightarrow 0$. Consequently $$\lim _{i
\rightarrow \infty} \delta _i \ln u_i(x_i) =0$$ which implies $\lim _{i
\rightarrow \infty}  u_i(x_i)^{ \delta _i} =1$.
\stop

\section{Ruling out Case I}
\label{section:case1}

In section \ref{section:iso}, we reduced the possible blow-up phenomenon of
$\{u_i\}$ which are solutions of equation (\ref{eq:main}) into two cases. In
this section we are going to show that case I can not happen, in
the next section we will rule out case II and hence complete the
proof of Theorem \ref{thm:main}.

\noindent
{\bf Case I:} 
The sequence $\{u_i\}$ has isolated blow-up points $x_1, x_2, ... \in M$.

Suppose $x_1, x_2, ...$ are all simple blow-up
points. Choose $P \in M \setminus \{ x_1, x_2,...\}
$. On any compact
subset $\Omega$ of $M \setminus \{ x_1, x_2, ...\}  $ containing
$P$, since $x_1, x_2, ...$ are isolated blow-up points, $u_i$ is
bounded above by some constant independent of $i$, so on $\Omega$ the
standard Harnack inequality holds for $\{ u_i \}$. Then by Proposition
\ref{prop:simpleestimates} and the Harnack inequality $u_i(P) \rightarrow 0$. In addition,
the Harnack inequality also holds for
$\frac{u_i}{u_i(P)} $. In other words, for some constant $C$
independent of $i$,
$$\max _{\Omega} \frac{u_i}{u_i(P)} \leq C \min _{\Omega}
\frac{u_i}{u_i(P)} \leq C \frac{u_i(P)}{u_i(P)}=C. $$  
Since $u_i$ satisfies (\ref{eq:main}), 
$$\Delta _{g} \left (\frac{u_i}{u_i(P)} \right) + u_i(P)^{p_i-1} K_i
\left (\frac{u_i}{u_i(P)} \right )^{p_i}=0 .$$ 
By the standard elliptic estimates, $\frac{u_i}{u_i(P)}$ has uniform
$C^{2, \alpha}$-norm on $\Omega$. So on $\Omega$, $\frac{u_i}{u_i(P)} \rightarrow
H$ in $C^2$-norm where $H$ satisfies $\Delta _{g} H=0$. Since $\Omega$ is arbitrary, $H$ satisfies $\Delta _{g} H=0$ on $M
\setminus \{ x_1, x_2, ...\}  $. By the fact $u_i >0$ we know that
$H(x) \geq 0$, then the maximal principle gives $H >0$ on $M
\setminus \{ x_1, x_2, ...\} $. Thus by the removable
singularity theorems of harmonic functions $H$ is a constant.  

Since we assume $x_1$ is a simple blow-up point, there exists a
sequence of points $\{ x^{(i)} _1 \}$ approaching $x_1$ such that for
the coordinates $z$ centered at each $\{ x^{(i)} _1 \}$, the function 
$|z|^{\frac{2}{p_i-1}}\bar{u}_i(|z|)$ is strictly decreasing in $|z|$ for
$Ru_i(x_i)^{-\frac{p_i-1}{2}} \leq |z| \leq r_0$. In particular it is
decreasing for
$\frac{r_0}{2} \leq |z| \leq r_0$ when $i$ is sufficiently large. This
implies that $|z|^{\frac{2}{p_i-1}}\bar{u}_i(|z|) u_i(P)^{-1}$ is
strictly decreasing in $|z|$ for $\frac{r_0}{4} \leq |z| \leq r_0$. 
If $\frac{r_0}{2} \leq |z| \leq r_0$, then the corresponding point
$\exp_{x^{(i)} _1 } z$ is at least distance $\frac{r_0}{4}$ from $x_1$
because $\{ x^{(i)}_1 \}$ approaches $x_1$.
Thus $\bar{u}_i(|z|) u_i(P)^{-1}$ converges in $C^2$-norm to $H$, which
is a constant. Consequently $|z|^{\frac{2}{p_i-1}}\bar{u}_i(|z|)
u_i(P)^{-1}$ converges in $C^2$-norm to $|z|^{\frac{2}{p_i-1}}H$ which
is strictly increasing in $|z|$. This is a contradiction.

Therefore there must be a point in $\{ x_1, x_2, ...\}$ which
is not a simple blow-up point, without loss of generality we assume it
to be $x_1$. To simplify the notations we are going to rename it to
be $x_0$. Let $ x_i $ be the local maximum points of $u_i$ such that
$\lim _{i \rightarrow \infty}x_i=x_0$. Let $z$ be the local coordinate
system centered at each $x_i$. Since $x_0$ is not a simple blow-up point, as
a function of $|z|$, $|z|^{\frac{2}{p_i-1}}\bar{u}_i(|z|)$ has a
second critical point at $|z|=r_i$ where $r_i \rightarrow 0$. Let
$y=\frac{z}{r_i}$ and define
$v_i(y)=r_i^{\frac{2}{p_i-1}}u_i(r_iy)$. Then $v_i(y)$ satisfies
\begin{equation}
\label{eq:rescaled-case1}
\Delta _{g^{(i)}} v_i(y)+\widetilde{K}_i(y)v_i(y)^{p_i}=0 
\end{equation} 
where $g^{(i)}(y) =g_{\alpha \beta}(r_iy)dy^{\alpha}dy^{\beta}$ and
$\widetilde{K}_i(y)=K_i(r_iy)$.\\ 

\noindent
By this definition 
$|y|=1$ is the second critical point of $|y|^{\frac{2}{p_i-1}}\bar{v}_i(|y|)$. 
By Proposition \ref{prop:iso}, for $0 \neq |z| \leq \sigma$, $u_i(z) \leq C
|z|^{-\frac{2}{p_i-1}}$ where $\sigma$ is a positive constant. Then
since $\frac{\sigma}{r_i} \rightarrow \infty $, $|v_i(y)| \leq C |y|^{-\frac{2}{p_i-1}}$ for $|y| \neq
0$. Therefore by the same argument as before we know that $v_i(y)$
converges in $C^2$-norm on $\mathbf{R}^n \setminus \{0\}$ to some function
$v$ which satisfies  $\Delta v + K(x_0)v^{\frac{n+2}{n-2}}=0$ where
here and in the rest of the proof $\Delta$ is the Euclidean Laplacian and $K$ is the limit function of
$\{K_i\}$.

\noindent
If $0$ is not a blow-up point of $\{v_i\}$, then $v$
satisfies $\Delta v + K(x_0)v^{\frac{n+2}{n-2}}=0$ on
$\mathbf{R}^n$. Since $r_i > Ru_i(x_i)^{-\frac{p_i-1}{2}}$, $v_i(0)
>R^{\frac{2}{p_i-1}}$ and hence $v>0$. This implies that $v$ is the
standard spherical solution and $|y|^{\frac{n-2}{2}}v(|y|)$ only has
one critical point.
On the other hand, $|y|^{\frac{2}{p_i-1}}\bar{v}_i(|y|)$ has two
critical points $$|y|=1 \hspace{.1in} \text{ and } \hspace{.1in} |y_i|=\left (
\frac{2n(n-2)}{(np_i-n-2p_i)K_i(x_i)} \right )^{\frac{1}{2}} u_i(x_i)^{-\frac{p_i-1}{2}}r_i^{-1}. $$
If $|y_i| \rightarrow 0$, then letting $i\rightarrow \infty$ we have
\begin{eqnarray*}
 0 & = & \big (|y|^{\frac{n-2}{2}}v(|y|) \big )(0) \,\, = \,\, \lim _{i \rightarrow
\infty} \big ( |y|^{\frac{2}{p_i-1}}\bar{v}_i(|y|) \big )(|y_i|) \\
& = & 2^{-\frac{n-2}{2}}\left (
\frac{n(n-2)}{K(x_0)} \right )^{\frac{n-2}{4}},
\end{eqnarray*}
this is a contradiction. So $|y_i|$ doesn't converge to $0$, then it
implies that $ |y|^{\frac{n-2}{2}}v(|y|)$ should also have two
critical points, which is a contradiction too.
Thus
$0$ is a blow-up point for $\{v_i\}$ and furthermore by the
construction a simple blow-up point. 

Choose a point $\bar{y}$ with $|\bar{y}|=1$. On any compact
subset $\Omega$ of $\mathbf{R}^n \setminus \{0\}$ which contains $\bar{y}$, similar to the
proof of Lemma \ref{lemma:simpleharnack} we have $\max _{\Omega} v_i \leq C(\Omega) \min
_{\Omega} v_i$, so $$\max _{\Omega} \frac{v_i}{v_i(\bar{y})} \leq C(\Omega) \min
_{\Omega} \frac{v_i}{v_i(\bar{y})} \leq C(\Omega).$$
Thus we can conclude that $\frac{v_i}{v_i(\bar{y})}$ converges in
$C^2$-norm on $\Omega$ to a function $h$. Additionally, since
$v_i(\bar{y}) \rightarrow 0$ by Proposition \ref{prop:simpleestimates} and
$$  \Delta _{g^{(i)}} \left(\frac{v_i}{v_i(\bar{y})} \right) +
v_i(\bar{y})^{p_i-1}\widetilde{K}_i(y) \left (\frac{v_i}{v_i(\bar{y})}
\right )  ^{p_i}=0,       $$
$h$ satisfies $\Delta h =0$ on $\Omega$. Therefore
$\Delta h =0$ on $\mathbf{R}^n \setminus \{0\}$ since $\Omega$ is
arbitrary.

\noindent
Because $0$ is a simple blow-up point for $\{v_i\}$, it is a
non-removable singularity for $h$. So the singular part of $h$ has the form $b|y|^{2-n}$
for some constant $b$. Now since $h-b|y|^{2-n}$ is harmonic on
$\mathbf{R}^n$, it is a constant. So we can write $h(y)=a+b|y|^{2-n}$
for constants $a$ and $b$. Since $h(\bar{y})=\lim _{i \rightarrow
\infty} \frac{v_i(\bar{y})}{v_i(\bar{y})}=1$ and $|\bar{y}|=1$, we
have $a+b=1$. Because $|y|=1$ is a critical point of
$|y|^{\frac{2}{p_i-1}}\frac {\bar{v}_i(|y|)}{v_i(\bar{y})}$, it is
also a critical point of
$|y|^{\frac{n-2}{2}}h(|y|)=a|y|^{\frac{n-2}{2}}+b|y|^{-\frac{n-2}{2}}$.
Taking the derivative at $|y|=1$ we have $a \frac{n-2}{2} - b
\frac{n-2}{2}=0$. So $a=b=1$ and $h=\frac{1}{2} + \frac{1}{2}|y|^{2-n}$.

Next we are going to apply the Pohozaev identity (\ref{eq:pohozaev}) to equation
(\ref{eq:rescaled-case1}). Since $g$ is locally conformally flat, we can write
$g(z)=\lambda ^{\frac{4}{n-2}} (z) dz^2$. Hence
we can write $g^{(i)}(y)= \lambda ^{\frac {4}{n-2}} (r _i y) dy^2$. We are going to use  $\lambda _i (y)$ to denote $ \lambda (r _i y)$. Let $X=\sum _{j} y^j \frac{\partial}{\partial y^j}$, the
Pohozaev identity becomes
\begin{equation}
\label{eq:pohozaev-case1}
\frac{n-2}{2n} \int _{B_{\sigma}} X(R_i) dv _{g_i} = \int _{\partial
B_{\sigma}} T_i(X, \nu _i) d \Sigma _i
\end{equation}
where 
\begin{eqnarray*}
g_i(y) & = & v_i(y)^{\frac{4}{n-2}}g^{(i)}(y) \,\,  = \,\,  (\lambda
_i v_i)^{\frac{4}{n-2}}d^2y, \\
R_i(y) & = & R (g_i)  \,\, = \,\,  c(n)^{-1} \widetilde{K} _i v_i^{-\delta _i},\\
dv_{g_i} & = & v_i(y)^{\frac{2n}{n-2}}dv_{g^{(i)}} \,\,  = \,\,  (\lambda _i v_i)^{\frac{2n}{n-2}}dy, \\  
\nu _i & = & (\lambda _i v_i)^{-\frac{2}{n-2}}\sigma ^{-1}\sum _{j}
y^j \frac{\partial}{\partial y^j} \\
& & \text{is the unit
outer normal vector on } \partial B_{\sigma} \text{ with respect to }
g_i,\\
d \Sigma _i & = & (\lambda _i v _i)^{\frac{2(n-1)}{n-2}}d
\Sigma _{\sigma} \\
& &  \text{where } d \Sigma _{\sigma} \text{ is the
surface element of the standard } S^{n-1}(\sigma),\\
T_i & = & \Ric (g_i)-n^{-1}R(g_i)g_i.
\end{eqnarray*}

We divide both sides of (\ref{eq:pohozaev-case1}) by $v_i^2(\bar{y})$. The right hand side
becomes 

\begin{eqnarray}
\label{eq:case1-Tterm}
 & & \frac{1}{v_i^2(\bar{y})}\int _{\partial
B_{\sigma}} T_i(X, \nu _i) d \Sigma _i \nonumber \\
& = & \frac{1}{v_i^2(\bar{y})}\int _{\partial
B_{\sigma}}\left ( \Ric (g_i)-n^{-1}R(g_i)g_i \right )(X, 
\nu _i) d \Sigma _i \nonumber  \\
& = & \frac{1}{v_i^2(\bar{y})}\int _{\partial
B_{\sigma}} 
\bigg ( \Ric 
  \left ( 
    \left (\lambda _i v_i \right )^{\frac{4}{n-2}}g_0 
  \right ) \nonumber \\
& & - n^{-1} R 
  \left( \left(\lambda _i v_i \right)^{\frac{4}{n-2}}g_0 \right) 
  \left( \lambda _i v_i \right )^{\frac{4}{n-2}}g_0 
\bigg )(X, \nu_0)
(\lambda _i v_i)^2   d \Sigma _{\sigma} \nonumber \\
\nonumber \\
& = &  \int _{\partial
B_{\sigma}} 
\left( \frac{\lambda _i v_i}{v_i(\bar{y})}\right)^2 
\Bigg [ \Ric
  \left ( 
     \left(\frac{\lambda _i v_i}{v_i(\bar{y})}\right)^{\frac{4}{n-2}} g_0
  \right) \\
  & & - n^{-1} R 
  \left ( 
    \left(\frac{\lambda _i
v_i}{v_i(\bar{y})}\right)^{\frac{4}{n-2}}g_0 
  \right) 
  \left( \frac{\lambda _i v_i}{v_i(\bar{y})} \right ) 
  ^{\frac{4}{n-2}}g_0 
\Bigg ] (X, \nu _0)    d \Sigma_{\sigma} \nonumber 
\end{eqnarray} 
 where  $g_0$  denotes the Euclidean metric and $\nu _0 =\sigma ^{-1}\sum _{j}
y^j \frac{\partial}{\partial y^j}$ is the unit outer normal on
$\partial B_{\sigma}$ with respect to 
the Euclidean metric $ g_0$. 

\noindent
When $i \rightarrow \infty$, for $|y|=\sigma$, $\lambda _i(y) =
\lambda (r _i y) \rightarrow \lambda (x_0)$, without loss of
generality we can assume it to be $1$. Thus when $i$ goes to $\infty$,
(\ref{eq:case1-Tterm}) converges to 
\begin{eqnarray}
\label{eq:case1-Tterm-followup}
& & \int _{\partial
B_{\sigma}} h^2  \bigg (\Ric \left (h^{\frac{4}{n-2}}g_0 \right )  - n^{-1}
R \left (h^{\frac{4}{n-2}}g_0 \right )h^{\frac{4}{n-2}}  g_0 \bigg )(X, \nu _0)    d \Sigma _{\sigma}       \nonumber \\
& = &  \int _{\partial
B_{\sigma}} h^2 \cdot (n-2) h^{\frac{2}{n-2}} \bigg [\Hess \left
(h^{-\frac{2}{n-2}} \right) (X, \nu _0)\nonumber \\
& & -\frac{1}{n}\Delta \left (h^{-\frac{2}{n-2}} \right )g_0 (X, \nu _0)\bigg ]   \,\,  d \Sigma _{\sigma}
\nonumber \\
& = & (n-2)\sigma^{-1} \int _{\partial
B_{\sigma}} h^{\frac{2(n-1)}{n-2}} \cdot \\
& & \bigg [ \Hess \left
(h^{-\frac{2}{n-2}} \right )(X,X)  -\frac{1}{n}\Delta \left (h^{-\frac{2}{n-2}} \right )\sigma ^2  \bigg ] \,\,  d \Sigma _{\sigma} \nonumber
\end{eqnarray}

\noindent
By the expression of $h$
\begin{equation*}
h^{-\frac{2}{n-2}} \,\, = \,\, \left( \frac{1}{2} (1+|y|^{2-n}) \right )^{-\frac{2}{n-2}}
\,\, = \,\, 2^{\frac{2}{n-2}}|y|^2 - \frac{2^{\frac{n}{n-2}}}{n-2}|y|^{n} +
O \left (|y|^{2(n-1)} \right )
\end{equation*}


\noindent
Then by direct computation 
$$\Hess \left( 2^{\frac{2}{n-2}}|y|^2 - \frac {2^{\frac{n}{n-2}}}{n-2}
|y|^{n}  \right )(X,X) -
\frac{1}{n}\Delta \left ( 2^{\frac{2}{n-2}}|y|^2 - \frac {2^{\frac{n}{n-2}}}{n-2}
|y|^{n} \right ) \sigma ^2 $$
$= -2^{\frac{n}{n-2}}(n-1)  \sigma ^n $

\noindent
Therefore $$ \Hess \left (h^{-\frac{2}{n-2}} \right )(X,X)
-\frac{1}{n}\Delta \left (h^{-\frac{2}{n-2}} \right )\sigma ^2 =  -2^{\frac{n}{n-2}}(n-1) \sigma
^n + O \left (\sigma ^{2(n-1)} \right )  . $$    
Also we know
 \begin{eqnarray*}
 h^{\frac{2(n-1)}{n-2}}
& = & \left (\frac
{1}{2} \right )^{\frac{2(n-1)}{n-2}}|y|^{-2(n-1)}\left (1+ O(|y|^{n-2}) \right ).
\end{eqnarray*}

\noindent
So we can conclude that (\ref{eq:case1-Tterm-followup}) 
\begin{eqnarray}
\label{eq:case1-Tterm-<0}
& = & -\frac{1}{2}(n-1)(n-2) \sigma^{-1} \int _{\partial B_{\sigma} }
\left (|y|^{-2(n-1)}+O(|y|^{-n})\right ) \cdot \nonumber \\
& & \left ( |y|^{n}+O(|y|^{2(n-1)})
\right )\sigma^{n-1}d\Sigma
_1 \nonumber \\
& = &  -\frac{1}{2}(n-1)(n-2) + O(\sigma ^{n-2}) \nonumber \\
& < & 0 
\end{eqnarray}
when we choose $\sigma$ to be sufficiently small. 

On the other hand, after being divided by $v_i^2(\bar{y})$, the left
hand side of (\ref{eq:pohozaev-case1}) is 
$$
\frac{n-2}{2n}c(n)^{-1} \frac{1}{v^2_i(\bar{y})}\int _{B_{\sigma}}
X(\widetilde{K}_iv_i^{-\delta _i})(\lambda _i v_i)^{\frac{2n}{n-2}}dy.
$$
We write
\begin{eqnarray}
\label{eq:case1-X(R)term}
& & \frac{1}{v^2_i(\bar{y})}\int _{B_{\sigma}}
X(\widetilde{K}_iv_i^{-\delta _i })(\lambda _i
v_i)^{\frac{2n}{n-2}}dy \nonumber \\
 & =
& \frac{1}{v^2_i(\bar{y})}\int _{B_{\sigma}}
X(\widetilde{K}_i)v_i^{p_i +1 }\lambda _i ^{\frac{2n}{n-2}}dy
-\frac{\delta _i}{v^2_i(\bar{y})}
  \int _{B_{\sigma}} \widetilde{K}_i\lambda _i ^{\frac{2n}{n-2}}v_i^{p_i}X(v_i)
dy .
\end{eqnarray}
The second term 

\begin{eqnarray*}
& = & -\frac{\delta _i}{p_i+1}\frac{1}{v^2_i(\bar{y})}
   \int _{B_{\sigma}} \widetilde{K}_i\lambda _i ^{\frac{2n}{n-2}}X(v_i^{p_i+1})
dy \\
& = &  -\frac{\delta _i}{p_i+1}\frac{1}{v^2_i(\bar{y})} \int
_{B_{\sigma}} \left (\diver(\widetilde{K}_i \lambda _i ^{\frac{2n}{n-2}}v_i^{p_i+1} X) -
\widetilde{K}_i \lambda _i ^{\frac{2n}{n-2}}v_i^{p_i+1} \diver X \right ) dy \\
& &  + \frac{\delta
_i}{p_i+1}\frac{1}{v^2_i(\bar{y})}  \int _{B_{\sigma}} \lambda _i
^{\frac{2n}{n-2}}v_i^{p_i+1}X(\widetilde{K}_i) dy \\
& & +  \frac{\delta
_i}{p_i+1}\frac{1}{v^2_i(\bar{y})}  \int _{B_{\sigma}} \widetilde{K}_i
v_i^{p_i+1} X (\lambda _i ^{\frac{2n}{n-2}} )
dy \\
& = &  -\frac{\delta _i}{p_i+1}\frac{\sigma}{v^2_i(\bar{y})} \int
_{\partial B_{\sigma}} \widetilde{K}_i \lambda _i
^{\frac{2n}{n-2}}v_i^{p_i+1} d \Sigma _{\sigma}  + \frac{\delta
_i}{p_i+1}\frac{1}{v^2_i(\bar{y})} \cdot \\
& & \int _{B_{\sigma}}  \widetilde{K}_i \lambda _i
^{\frac{2n}{n-2}}v_i^{p_i+1} \left (n + X(\ln \widetilde{K}_i)+ \frac{2n}{n-2}
X(\ln \lambda _i) \right ) dy
\end{eqnarray*}

\noindent
Since $X=r\frac{\partial}{\partial r}$ and $\frac{\partial}{\partial
r}(\ln \widetilde{K}_i)$, $\frac{\partial}{\partial
r}(\ln \lambda _i)$ are uniformly bounded, we can choose $\sigma$ to be
small (independent of $i$) to make $n + X(\ln \widetilde{K}_i )+ \frac{2n}{n-2}
X(\ln \lambda _i)>0$.
\noindent
Because on $\partial B_{\sigma}$, $\frac{v_i}{v_i(\bar{y})}
\rightarrow h(\sigma)>0$ and $v_i \rightarrow 0$ uniformly,
$$
\frac{1}{v^2_i(\bar{y})} \int
_{\partial B_{\sigma}} \widetilde{K}_i \lambda _i
^{\frac{2n}{n-2}}v_i^{p_i+1} d \Sigma _{\sigma}   =   \int
_{\partial B_{\sigma}} \widetilde{K}_i \lambda _i
^{\frac{2n}{n-2}} \left ( \frac {v_i}{v_i(\bar{y})} \right )^2  v_i^{p_i-1} d \Sigma
_{\sigma}  \rightarrow  0 .
$$
Thus when $i \rightarrow \infty$, the limit of the second term of (\ref{eq:case1-X(R)term}) is greater than or equal
to $0$.  

As will be proved in Proposition \ref{prop:X(K)term}, when the dimension $n=3, 4$, the limit of the first
term of (\ref{eq:case1-X(R)term})
$$\lim _{i \rightarrow \infty} \frac{1}{v^2_i(\bar{y})}\int _{B_{\sigma}}
X(\widetilde{K}_i)v_i^{p_i+1}\lambda _i^{\frac{2n}{n-2}}dy
=0.$$
This then implies that the limit of the left hand side of (\ref{eq:pohozaev-case1}) is greater than
or equal to $0$, which contradicts (\ref{eq:case1-Tterm-<0}). So we can rule
out Case I.

\begin{prop}
\label{prop:X(K)term}
When $n=3, 4$, $$\lim _{i \rightarrow \infty} \frac{1}{v^2_i(\bar{y})}\int _{B_{\sigma}}
X(\widetilde{K}_i)v_i^{p_i+1}\lambda _i^{\frac{2n}{n-2}}dy
=0.$$
\end{prop}

Before we prove Proposition \ref{prop:X(K)term}, we first need to carefully investigate the behaviour of $\widetilde{K}_i$.

By Proposition \ref{prop:simpleestimates} we have the following estimates:
\begin{itemize}
\item if \,\, $0
\leq |y| \leq 1$, \hspace{.1in} $v_i(y) \geq C v_i(0)\left(1+ \frac{\widetilde{K}_i(0)}{n(n-2)}v_i(0)^{\frac{4}{n-2}}|y|^2 \right )^{-\frac{n-2}{2}}$ 
\item  if \,\, $0 \leq |y| \leq Rv_i(0)^{-\frac{p_i-1}{2}}$,
\hspace{.1in} then\\ $v_i(y) \leq C v_i(0) \left (1+ \frac{\widetilde{K}_i(0)}{n(n-2)}v_i(0)^{p_i-1}|y|^2 \right )^{-\frac{n-2}{2}}$
\item if \,\, $\frac{R}{v_i(0)^{\frac{p_i-1}{2}}} \leq |y| \leq 1$,
\hspace{.1in} then \hspace {.1in}$v_i(y) \leq Cv_i(0)^{t _i}|y|^{-l_i}$
\hspace{.1in} where $l_i$, $t_i$ are chosen such that 
$l_i \rightarrow \frac{6(n-2)}{7}$,
and $t _i =1-\frac{(p_i-1)l_i}{2}$.
\end{itemize}

\begin{lemma}
\label{lemma:derivativeofK}
For any $j=1,2,...,n$, 
$$\bigg |\frac{\partial \widetilde{K}_i}{\partial y^j}(0) \bigg | \leq C
\left ( r_iv_i(0)^{-\frac{2}{n-2}+\frac{n-1}{2}\delta_i}+
v_i(0)^{2t_i} \right ) $$ 
\end{lemma}

\pf
Choose the conformal Killing vector field to be $X=\frac
{\partial}{\partial y^1}$, we have the Pohozaev identity
\begin{equation}
\label{eq:pohozaev-for-derivativeofK}
\frac{n-2}{2n} \int _{B_{\sigma}} X(R_i) dv_{g_i} = \int _{\partial B_{\sigma}} T_i(X, \nu _i) d \Sigma _i
\end{equation}
where
\begin{eqnarray*}
g_i(y) & = & v_i(y)^{\frac{4}{n-2}}g^{(i)}(y) \,\,  = \,\,    (\lambda
_i v_i)^{\frac{4}{n-2}}g_0 \\
& & \text{where } g_0 \text{ is the Euclidean metric,} \\
R_i(y) & = & R (g_i) \,\,  = \,\,  c(n)^{-1} \widetilde{K} _i v_i^{-\delta _i},\\
dv_{g_i} & = & v_i(y)^{\frac{2n}{n-2}}dv_{g^{(i)}} \,\,  = \,\,
(\lambda _i v_i)^
{\frac{2n}{n-2}}dy, \\  
\nu _i & = & (\lambda _i v_i)^{-\frac{2}{n-2}}\sigma ^{-1}\sum _{j}
y^j \frac{\partial}{\partial y^j} \\
& &  \text{is the unit
outer normal vector on } \partial B_{\sigma} \text{ with respect to } g_i,\\
d \Sigma _i & = & (\lambda _i v _i)^{\frac{2(n-1)}{n-2}}d \Sigma _{\sigma} \\
& & \text{where } d \Sigma _{\sigma}  \text{ is the
surface element of the standard } S^{n-1}(\sigma),\\
T_i & = & (n-2) (\lambda_i v_i)^{\frac{2}{n-2}} \left ( \Hess \big ( (\lambda _i
v_i)^{-\frac{2}{n-2}} \big)-\frac{1}{n}\Delta \big ((\lambda_i
v_i)^{-\frac{2}{n-2}} \big) g_0 \right ).
\end{eqnarray*}
Here $\Hess$ and $\Delta$ are taken with respect to the Euclidean
metric $g_0$ .

The left hand side of (\ref{eq:pohozaev-for-derivativeofK}) is 
\begin{eqnarray}
\label{eq:derivativeofK-lhs}
& &  \frac{n-2}{2n} \int _{B_{\sigma}} \frac{\partial}{\partial
y^1}(R_i) dv_{g_i} \nonumber \\
& = & \frac{n-2}{2n} c(n)^{-1} \int _{B_{\sigma}} \frac{\partial}{\partial
y^1}(\widetilde{K}_iv_i^{-\delta _i} )(\lambda_i v_i)^{\frac{2n}{n-2}}dy \nonumber \\
& = & \frac{n-2}{2n} c(n)^{-1} \int _{B_{\sigma}} \left (1+
\frac{\delta_i}{p_i+1} \right )\lambda_i^
{\frac{2n}{n-2}}v_i^{p_i+1}\frac{\partial \widetilde{K}_i}{\partial
y^1} dy \nonumber \\
& & +\frac{n-2}{2n}c(n)^{-1} \int _{B_{\sigma}} \frac{\delta_i}{p_i+1}
\widetilde{K}_i v_i^{p_i+1} \frac {\partial \lambda _i^
{\frac{2n}{n-2}}}{\partial y^1} dy  \nonumber\\
& & - \frac{n-2}{2n} c(n)^{-1}
\frac{\delta_i}{p_i+1}\int _{\partial B_{\sigma}}
\lambda_i^{\frac{2n}{n-2}}\widetilde{K}_iv_i^{p_i+1} \frac{y^1}{\sigma}
d \Sigma _{\sigma}
\end{eqnarray}
\noindent
By Proposition \ref{prop:simpleestimates}, the third term in (\ref{eq:derivativeofK-lhs}) is
bounded above by $$C\delta _i \cdot v_i(0)^{t_i(p_i+1)} \leq C\delta _iv_i(0)^{2t _i}
$$ since $t_i < 0$ and $v_i(0) \rightarrow \infty$.

\noindent
Same as in the proof of Proposition \ref{prop:u^delta}, the second term in (\ref{eq:derivativeofK-lhs}) is bounded above by \\
\begin{eqnarray*}
& & C \delta _i r_i\int _{|y| \leq \sigma} v_i(y)^{p_i+1} dy \\
& \leq & C \delta _i r_i \Bigg (\int _{|z| \leq Rv_i(0)^{-\frac{p_i-1}{2}}}\,\,
v_i(0)^{p_i+1}\,\,  dy\\
& & \hspace{.3in}  +\int _{ Rv_i(0)^{-\frac{p_i-1}{2}} \leq
|y| \leq \sigma} \left (v_i(0)^{t _i}|y|^{-l_i} \right )^{p_i+1} dy \Bigg )\\
& \leq & C \delta _i r_i \left ( v_i(0)^{p_i+1 - \frac{n}{2} (p_i-1)}
+ v_i(0)^{t_i(p_i+1)} \cdot v_i(0)^{-\frac{p_i-1}{2} \left (
n-l_i(p_i+1) \right )} \right)\\
& = & C \delta _i r_i v_i(0)^{p_i+1-\frac{p_i-1}{2}n} \hspace{.5in} (\text{ since }t _i + \frac{(p_i-1)l_i}{2}=1 )\\
& = & C \delta _i r_i v_i(0)^{(\frac{n}{2}-1)\delta_i}\\
& \leq & C \delta _i r_i \hspace {1.5in} \text{(by Proposition \ref{prop:u^delta})}
\end{eqnarray*}

\noindent
By the estimates almost identical to those of the right hand side of
(\ref{eq:pohozaev-X(R)}) we know that the right hand of (\ref{eq:pohozaev-for-derivativeofK}) decays in the rate of
$ v_i(0)^{2 t _i}$.

\noindent
Therefore the first term in (\ref{eq:derivativeofK-lhs}) which is $$\frac{n-2}{2n} c(n)^{-1} \int _{B_{\sigma}} \left(1+
\frac{\delta_i}{p_i+1}\right)\lambda_i^
{\frac{2n}{n-2}}v_i^{p_i+1}\frac{\partial \widetilde{K}_i}{\partial
y^1} dy$$ is bounded above by 
$
C(\delta _i v_i(0)^{2t_i}+ \delta _i r_i
+ 
v_i(0)^{2 t _i}) \leq C\left( \delta _i r_i
+ 
v_i(0)^{2 t _i}  \right )$.

\noindent
By the Taylor expansion
$$\frac{\partial \widetilde{K}_i}{\partial y^1}(y)= \frac{\partial \widetilde{K}_i}{\partial
y^1}(0) + \triangledown \left(\frac{\partial \widetilde{K}_i}{\partial
y^1} \right )(v) \cdot y
\hspace{.3in} \text{ for some } |v| \leq |y|.$$

\noindent
As in the proof of Proposition \ref{prop:u^delta} and also using the fact that
$\widetilde{K} _i(y)=K_i(r_iy)$,

\begin{eqnarray*}
& & \frac{n-2}{2n} c(n)^{-1} \int _{B_{\sigma}} \left (1+
\frac{\delta_i}{p_i+1} \right )\lambda_i^
{\frac{2n}{n-2}}v_i^{p_i+1} \Bigg | \triangledown \left (\frac{\partial
\widetilde{K}_i}{\partial y^1} \right )(v) \cdot y \Bigg | dy \\ & \leq & C r_i \int _{B_{\sigma}}
v_i^{p_i+1}|y|dy \\
& \leq & C r_i v_i(0)
^{-\frac{2}{n-2} +\frac{n-1}{2}\delta_i}  
\end{eqnarray*}
where the last inequality is proved in the same way as (\ref{eq:int-|z|u_i^p_i+1}).
 
\noindent
Thus we know 
\begin{eqnarray*}
& & \bigg |\frac{\partial \widetilde{K}_i}{\partial y^1}(0) \bigg |\int  _{B_{\sigma}} \lambda_i^
{\frac{2n}{n-2}}v_i^{p_i+1} dy\\
 & \leq & C \left (r_iv_i(0)^{-\frac{2}{n-2}+\frac{n-1}{2}\delta_i} + \delta_i r_i
+ 
v_i(0)^{2 t _i} \right )\\
&  \leq & C \left (r_iv_i(0)^{-\frac{2}{n-2}+\frac{n-1}{2}\delta_i} + r_i v_i(0)^{2t_i}
+ v_i(0)^{2t_i} \right )\\ 
& & \text{ (by inequality (\ref{eq:upper-on-delta}))}\\
& \leq & C \left (r_iv_i(0)^{-\frac{2}{n-2}+\frac{n-1}{2}\delta_i}+
v_i(0)^{2t_i} \right )
\end{eqnarray*} 
Then by (\ref{eq:lower-int-u^p+1}) 
$$
\bigg |\frac{\partial \widetilde{K}_i}{\partial y^1}(0) \bigg | \leq C
\left (r_iv_i(0)^{-\frac{2}{n-2}+\frac{n-1}{2}\delta_i} + v_i(0)^{2 t _i}
\right ).$$

The same estimate holds for $\big |\frac{\partial \widetilde{K}_i}{\partial
y^j}(0) \big |, \,\, j=2,...n $ as well.

\stop

Now we can prove Proposition \ref{prop:X(K)term}.

\pf 
By the estimates of $v_i$ as stated between Proposition \ref{prop:X(K)term} and Lemma
\ref{lemma:derivativeofK}, it is equivalent to proving 
 $$\lim _{i \rightarrow \infty} v^2_i(0)\int _{B_{\sigma}}
X(\widetilde{K}_i)v_i^{p _i+1 }\lambda _i^{\frac{2n}{n-2}}dy
=0.$$

When $n=3,4$,
\begin{eqnarray*}
 X(\widetilde{K}_i)(y) 
 & = & \left (\sum _{j} y^j \frac{\partial
\widetilde{K}_i}{\partial y^j} \right )(y) \\
& = & \left (\sum _{j} y^j \frac{\partial
\widetilde{K}_i}{\partial y^j} \right )(0) + \sum _{k} 
\frac{\partial}{\partial y^k} \left (\sum _{j}y^j \frac{\partial
\widetilde{K}_i}{\partial y^j} \right ) (0)y^k \\
 & & 
+ \sum _{k,l}
\frac{\partial ^2}{\partial y^k \partial y^l} \left (\sum _{j}y^j \frac{\partial
\widetilde{K}_i}{\partial y^j} \right )(0)y^ky^l+O(|y|^3)  \\
& = & \sum _{j} \frac{\partial
\widetilde{K}_i}{\partial y^j}(0) y^j + \sum _{j,k} \frac{\partial^2
\widetilde{K}_i}{\partial y^j \partial y^k}(0) y^jy^k + O(|y|^3) 
\end{eqnarray*}

\noindent
Since for $\kappa = 1$ or $\kappa =2$,
$$n-l_i(p_i+1)+\kappa \rightarrow -\frac{5n}{7}+\kappa <0,  $$
by similar calculation as in the proof of Proposition \ref{prop:u^delta} We have 
\begin{eqnarray}
& & \int _{B_{\sigma}} v_i^{p
_i+1 }|y|^{\kappa}dy \nonumber \\
& \leq & C \Bigg ( \int _  { |y| \leq Rv_i(0)^{-\frac{p_i-1}{2}}  }\,\,
                 v_i(0)^{p_i+1}   |y|^{\kappa} \,\, dy \nonumber \\
& &         +\int _{ Rv_i(0)^{-\frac{p_i-1}{2}} \leq |y| \leq \sigma} 
             \left  (v_i(0)^{t _i}|y|^{-l_i} \right )^{p_i+1} |y|^{\kappa} dy 
            \Bigg )
\nonumber \\
& \leq  &   Cv_i(0)^{p_i+1 -
\frac{(n+\kappa)(p_i-1)}{2}} 
\end{eqnarray}

\noindent
Then
\begin{eqnarray*}
 v^2_i(0) \bigg |\int _{B_{\sigma}} \frac{\partial^2
\widetilde{K}_i}{\partial y^j \partial y^k}(0) y^jy^k 
v_i^{p _i+1 }\lambda _i^{\frac{2n}{n-2}}dy \bigg |
& \leq & Cr_i^2v^2_i(0)\int _{B_{\sigma}}
v_i^{p _i+1 }|y|^2dy\\
& \leq & Cr_i^2v_i(0)^{\frac{2n-8}{n-2}+ \frac{n}{2}\delta _i }\\
& \rightarrow & 0 \hspace{.3in} \text{ as } i \rightarrow 0.
\end{eqnarray*}
By Proposition \ref{prop:u^delta} and Lemma \ref{lemma:derivativeofK},

\begin{eqnarray*}
& &  v^2_i(0) \bigg |\int _{B_{\sigma}} \frac{\partial
\widetilde{K}_i}{\partial y^j}(0) y^j
v_i^{p _i+1 }\lambda _i^{\frac{2n}{n-2}}dy \bigg | \\
& \leq & Cv^2_i(0) \left (r_iv_i(0)^{-\frac{2}{n-2}+\frac{n-1}{2}\delta
_i}+v_i(0)^{2t_i} \right ) \int _{B_{\sigma}}|y|v_i^{p_i+1}dy \\
& \leq & C \left (r_i v_i(0)^{2-\frac{4}{n-2}+(n-1)\delta _i
}+v_i(0)^{2-\frac{2}{n-2}+2t_i+\frac{n-1}{2}\delta_i} \right) \\
& \leq & C \left (r_iv_i(0)^{2-\frac{4}{n-2}}+v_i(0)^{2-\frac{2}{n-2}+2t_i} \right )
\end{eqnarray*}

\noindent
Since 
$t _i =1 -\frac{(p_i-1)l_i}{2} \rightarrow 1-\frac{2}{n-2}
\frac{6(n-2)}{7} = - \frac {5}{7}, $
$$
\lim _{i \rightarrow \infty}
\left (2-\frac{2}{n-2}+2t_i \right ) \,\,  = \,\,  2-\frac{2}{n-2} - \frac {10}{7}
\,\, < \,\, 0.
$$

\noindent
We also have $2-\frac{4}{n-2} \leq 0$.
Therefore 
\begin{eqnarray*}
& &  v^2_i(0) \bigg |\int _{B_{\sigma}} \frac{\partial
\widetilde{K}_i}{\partial y^j}(0) y^j
v_i^{p _i+1 }\lambda _i^{\frac{2n}{n-2}}dy \bigg |\\
 &  \leq &   C \left
(r_iv_i(0)^{2-\frac{4}{n-2}}+v_i(0)^{2-\frac{2}{n-2}+2t_i} \right ) \\
&  \rightarrow & 0
 \hspace{.3in} \text { as } \,\, i \rightarrow 0.
\end{eqnarray*}

\noindent
Lastly, 
\begin{eqnarray*}
v^2_i(0)  \int _{B_{\sigma}} v_i^{p
 _i+1 }|y|^3 dy  &= & v^2_i(0) \Bigg ( \int _  {  |y| \leq R v_i(0)^{ -  \frac{p_i-1}{2}  } } 
\,\, v_i(y)^{p_i+1} |y|^3 \,\,  dy  \\
& & + \int _{ Rv_i(0)^{-\frac{p_i-1}{2}} \leq
|y| \leq \sigma} v_i(y)^{p_i+1}|y|^3 \,\, dy \Bigg ) 
\end{eqnarray*}

\noindent
The first term
$$
v^2_i(0) \int _{|y| \leq \frac{R}{v_i(0)^{\frac{p_i-1}{2}}}}
v_i(y)^{p_i+1} |y|^3  dy \,\, \leq \,\, C v_i(0)^{2+p_i+1-\frac{(n+3)(p_i-1)}{2}}
\,\,  \rightarrow \,\, 0 $$
since $2+p_i+1-\frac{(n+3)(p_i-1)}{2} \rightarrow 2- \frac{6}{n-2} <0.$

\noindent
Because $\lim _{i \rightarrow \infty} \big (-l_i(p_i+1)+3+n \big ) = -\frac
{5n}{7}+3 >0$,
the second term 

\begin{eqnarray*}
& &  v^2_i(0)\int _{ \frac{R}{v_i(0)^{\frac{p_i-1}{2}}} \leq
|y| \leq \sigma} v_i(y)^{p_i+1}|y|^3 dy\\
 & \leq & C v^2_i(0)\int _{ \frac{R}{v_i(0)^{\frac{p_i-1}{2}}} \leq
                        |y| \leq \sigma} 
             \left  (v_i(0)^{t _i}|y|^{-l_i} \right )^{p_i+1} 
               |y|^3 dy \\
& \leq & C v_i(0)^{2+t_i(p_i+1)} \\
& \rightarrow & 0
\end{eqnarray*}
since $2+t_i(p_i+1) \rightarrow 2- \frac{5}{7} \frac{2n}{n-2} <0$.\\
Therefore $\,\, v^2_i(0)  \int _{B_{\sigma}} v_i^{p
 _i+1 }|y|^3 dy\,\, $ also converges to $0$. 

\noindent
Thus we have proved that $\lim _{i \rightarrow \infty} v^2_i(0)\int _{B_{\sigma}}
X(\widetilde{K}_i)v_i^{p _i+1 }\lambda _i^{\frac{2n}{n-2}}dy
=0 $ when $ n=3,4.$
\stop

\section{Ruling Out Case II}
\label{section:case2}

In section \ref{section:iso}, we have reduced Case II to the following:\\
There is a sequence of functions $\{v_i\}$, each satisfies  
$$
\Delta _{g^{(i)}} v_i + K_i(\sigma _i y)v_i^{p_i}=0
$$
where $g^{(i)}(y)=g_{\alpha \beta}(\sigma _i y)dy^{\alpha}dy^{\beta}$ converges to the
Euclidean metric on compact subset of $\mathbf{R}^n(y)$. The sequence
$\{ v_i \}$ has isolated blow-up point(s) $\{0,...\}$. 

If $0$ is not a simple blow-up point, then we can use the same
argument as in the previous section to rescale the function and get a
contradiction by examining both sides of the Pohozaev identity.

Thus $0$ must be a simple blow-up point for $\{v_i\}$, which satisfies
$\Delta _{g^{(i)}} v_i + K_i(\sigma _i y)v_i^{p_i}$\\$=0$. Then we can
study this sequence of $\{v_i\}$ in the same way as we studied the sequence of
solutions $\{v_i\}$ for equation (\ref{eq:rescaled-case1}). We can apply almost exactly
the same argument and get a contradiction. The only difference is
the expression of $h=\lim _{i \rightarrow \infty}
\frac{v_i(y)}{v_i(\bar{y})}$.

\noindent
In Case I, we know $h$ satisfies $\Delta h = 0$ on
$\mathbf{R}^n \setminus \{0\}$, $h(1)=1$, and from the construction has a second critical point at
$|y|=1$, which implies that $h=\frac{1}{2}(1+|y|^{2-n})$.  

\noindent
But here we only know $h$ satisfies $\Delta h = 0$ where it is regular and $h(1)=1$, but don't know
whether it has a second critical point. What we do know though is that
$0$ is not the only blow-up point of $\{v_i\}$. This is true because as
defined in section \ref{section:iso}, $v_i(y) = \sigma _i^{\frac{2}{p_i-1}}
u_i(\sigma _i y)$ and by the choice of $\sigma _i$, there exists $\{
y_{2,i} \}$ such that $|y_{2,i}|=1$ and  
$$v_i(y_{2,i}) = \sigma _i^{\frac{2}{p_i-1}}u_i(x_{2,i}) \geq
\left (\frac{R}{u_i(x_{2,i})^{\frac{p_i-1}{2}}} \right )^{\frac{2}{p_i-1}}u_i(x_{2,i})
= R^{\frac{2}{p_i-1}}.$$
Suppose $0$
is the only blow-up point for $\{v_i\}$, then the Harnack inequality
holds on any compact subset $\Omega$ of $\mathbf{R}^n(y)$ which contains $\partial
B_{\bar{r}}$ and a neighborhood of $y_2 = \lim _{i \rightarrow \infty}
y_{2,i}$, where $\bar{r}$ is chosen as in Proposition
\ref{prop:simpleestimates}. Therefore we know
$$ v_i(y_{2,i}) \leq \max _{\Omega} v_i \leq C \inf _{\Omega} v_i \leq
 C \inf _{\partial
B_{\bar{r}}} v_i \rightarrow 0 \hspace{.3in} \text {as } i \rightarrow
\infty
$$ by Proposition \ref{prop:simpleestimates}. This is a contradiction. 

\noindent
Thus $\{v_i\}$ has two or more 
blow-up points $\{0, y_2,...\}$, and
hence $h$ also has blow-up points $\{0, y_2,...\}$. Since $h$ is harmonic
everywhere else, we can write
$$h(y)=c_1|y|^{2-n}+ c_2 |y -y _2|^{2-n}  + \widetilde{h}(y)$$ where
$c_1, c_2>0$ are constants and $\widetilde{h}(y)$ is harmonic on
$\mathbf{R}^n \setminus \{y_3,...\}$ (if $h(y)$ has blow-up points
$y_3,...$ 
other than $0$ and $y_2$). By the harmonicity of $\widetilde{h}$,
the Harnack inequality and 
the
maximal principle, the infimum of $\widetilde{h}$ is
approached when $y$ goes off to $\infty$. Now since we know $h>0$ and
$\lim _{|y| \rightarrow \infty} c_1|y|^{2-n}= \lim _{|y| \rightarrow
\infty} c_2|y-y_2|^{2-n}=0$, the infimum of $\widetilde{h}$ must be
non-negative. Thus
when $|y|$ is small, $c_2 |y -y _2|^{2-n}
+ \widetilde{h}(y)  >0, $ i.e., 
near $0$,
$$h(y) = c_1|y|^{2-n}+  A + O(|y|),  \hspace {.3in} \text{ where }
\hspace{.1in} A>0.$$

\noindent
Then we can analyze $\frac{1}{v^2_i(\bar{y})}\int _{\partial
B_{\sigma}} T_i(X, \nu _i) d \Sigma _i$ in the same way as we did for
(\ref{eq:case1-Tterm}). The positive ``mass'' term $A >0 $ guarantees that $$\lim _{i
\rightarrow \infty} \frac{1}{v^2_i(\bar{y})}\int _{\partial
B_{\sigma}} T_i(X, \nu _i) d \Sigma _i <0.$$ By exactly the same
argument as in the previous section we can also show  
$$
\lim _{i \rightarrow \infty} \frac{n-2}{2n}c(n)^{-1} \frac{1}{v^2_i(\bar{y})}\int _{B_{\sigma}}
X(R_i)dv_{g_i}
\geq 0.
$$ This is a contradiction, so Case II is also ruled out.

Thus we have finished the proof of Theorem \ref{thm:main}.
 
\vspace{.2in}

When $R(g)>0$ and $K>0$, we can similarly define isolated blow-up points and
simple blow-up points for $\{u_i\}$ which satisfies $\Delta _{g} u_i -
c(n)
R(g)u_i +Ku_i^{p_i}=0$. After slight modification we have the same
estimates as those in Proposition \ref{prop:simpleestimates}. If the blow-up is not simple,
then it is either not isolated blow-up or it is isolated but not
simple blow-up.

\noindent
If the blow-up is isolated but not simple, we can rescale the
function and metric as in Case I of the scalar flat case to reduce
it to the simple blow-up case. Then a contradiction follows from the Pohozaev
identity as in the scalar flat case.

\noindent
If the blow-up is not isolated, we can first rescale the function in
the same way as in case II of the scalar flat
case to reduce it to the isolated blow-up cases. Then we can use
almost the 
identical argument as in the scalar flat case to rule it out.

\noindent
Thus the possible blow-up could only be simple. This completes the
proof of Theorem \ref{thm:cor}.

\section *{Acknowledgement}

This paper is based on part of the author's dissertation \cite{YY} at
Stanford University. I would like to thank my advisor Professor Richard Schoen for his help and support.  

\bibliographystyle{plain}
    \bibliography{thesis}

\begin{thebibliography}{10}

\bibitem{Au}
Thierry Aubin.
\newblock \'{E}quations diff\'erentielles non lin\'eaires et probl\`eme de
  {Y}amabe concernant la courbure scalaire.
\newblock {\em J. Math. Pures Appl. (9)}, 55(3):269--296, 1976.

\bibitem{CGS}
Luis~A. Caffarelli, Basilis Gidas, and Joel Spruck.
\newblock Asymptotic symmetry and local behavior of semilinear elliptic
  equations with critical {S}obolev growth.
\newblock {\em Comm. Pure Appl. Math.}, 42(3):271--297, 1989.

\bibitem{ES}
Jos{\'e}~F. Escobar and Richard~M. Schoen.
\newblock Conformal metrics with prescribed scalar curvature.
\newblock {\em Invent. Math.}, 86(2):243--254, 1986.

\bibitem{GT}
David Gilbarg and Neil~S. Trudinger.
\newblock {\em Elliptic partial differential equations of second order}.
\newblock Classics in Mathematics. Springer-Verlag, Berlin, 2001.
\newblock Reprint of the 1998 edition.

\bibitem{LZ}
Yanyan Li and Meijun Zhu.
\newblock Yamabe type equations on three-dimensional {R}iemannian manifolds.
\newblock {\em Commun. Contemp. Math.}, 1(1):1--50, 1999.

\bibitem{S1}
Richard Schoen.
\newblock Conformal deformation of a {R}iemannian metric to constant scalar
  curvature.
\newblock {\em J. Differential Geom.}, 20(2):479--495, 1984.

\bibitem{S2}
Richard Schoen.
\newblock Topics in differential geometry.
\newblock {\em Lecture Notes, Stanford University}, 1988.

\bibitem{S5}
Richard~M. Schoen.
\newblock The existence of weak solutions with prescribed singular behavior for
  a conformally invariant scalar equation.
\newblock {\em Comm. Pure Appl. Math.}, 41(3):317--392, 1988.

\bibitem{Tr}
Neil~S. Trudinger.
\newblock Remarks concerning the conformal deformation of {R}iemannian
  structures on compact manifolds.
\newblock {\em Ann. Scuola Norm. Sup. Pisa (3)}, 22:265--274, 1968.

\bibitem{YY}
Yu~Yan.
\newblock Some geometric problems involving conformal deformation of metrics.
\newblock {\em Dissertation, Stanford University}, 2004.

\end{thebibliography}

\end{document}